\documentclass[12pt]{article}
\usepackage{amssymb}

\let\ds=\displaystyle

\def\hybrid{\topmargin 0pt      \oddsidemargin 0pt
        \headheight 0pt \headsep 0pt
        \textwidth 17.5cm
        \textheight 25cm
        \voffset=-0.2cm
        \hoffset=-0.4cm
        \marginparwidth 0.0in
        \parskip 5pt plus 1pt   \jot = 1.5ex}
\catcode`\@=11
\def\marginnote#1{}

\newcount\hour
\newcount\minute
\newtoks\amorpm
\hour=\time\divide\hour by60
\minute=\time{\multiply\hour by60 \global\advance\minute by-\hour}
\edef\standardtime{{\ifnum\hour<12 \global\amorpm={am}%
        \else\global\amorpm={pm}\advance\hour by-12 \fi
        \ifnum\hour=0 \hour=12 \fi
        \number\hour:\ifnum\minute<10 0\fi\number\minute\the\amorpm}}
\edef\militarytime{\number\hour:\ifnum\minute<10 0\fi\number\minute}

\def\draftlabel#1{{\@bsphack\if@filesw {\let\thepage\relax
   \xdef\@gtempa{\write\@auxout{\string
      \newlabel{#1}{{\@currentlabel}{\thepage}}}}}\@gtempa
   \if@nobreak \ifvmode\nobreak\fi\fi\fi\@esphack}
        \gdef\@eqnlabel{#1}}
\def\@eqnlabel{}
\def\@vacuum{}
\def\draftmarginnote#1{\marginpar{\raggedright\scriptsize\tt#1}}

\def\draft{\oddsidemargin -0.1truein
        \def\@oddfoot{\sl preliminary draft\ {\tt\filename} \hfil
        \rm\thepage\hfil\sl\today\quad\militarytime}
        \let\@evenfoot\@oddfoot \overfullrule 3pt
        \let\label=\draftlabel
        \let\marginnote=\draftmarginnote
   \def\@eqnnum{{\rm (\theequation)}\rlap{\kern\marginparsep\tt\@eqnlabel}%
\global\let\@eqnlabel\@vacuum}  }
\def\numberbysection{\@addtoreset{equation}{section}
        \def\theequation{\thesection.\arabic{equation}}}
\numberbysection

\renewcommand{\theequation}{\thesection.\arabic{equation}}
\newcommand{\l@qq}[2]{\addvspace{2em}
 \hbox to\textwidth{\hspace{1em}\bf #1 \dotfill #2}}

\def\appname{Appendix}
\newcounter{app}
\def\theapp{\Alph{app}}
\def\app{\par
   \addvspace{4ex}
   \@afterindentfalse
  \secdef\@app\@dapp}
\def\@app[#1]#2{\ifnum \c@secnumdepth >\m@ne
        \refstepcounter{app}
        \addcontentsline{toc}{app}{\theapp
        \hspace{1em}#1}\else
      \addcontentsline{toc}{app}{ #1}\fi
   {\parindent \z@ \raggedright
    \large \bf \appname~\theapp .
   \large  \bf \hspace{1em}    #2}\nobreak
   \vskip -2ex   \noindent
\setcounter{equation}{0}
\def\theequation{\Alph{app}.\arabic{equation}}}
\def\@dapp#1{%
{\parindent \z@ \raggedright  \bf #1}\par\nobreak}
\def\l@app#1#2{\addpenalty{\@secpenalty}%
   \addvspace{1em plus\p@}%
   \begingroup
   \@tempdima 3em
     \parindent \z@ \rightskip \@pnumwidth
     \parfillskip -\@pnumwidth
     { \bf
     \leavevmode
     #1\hfil \hbox to\@pnumwidth{\hss #2}}\par
     \nobreak
   \endgroup}
\newcounter{sapp}[app]
\def\thesapp{\Alph{app}.\arabic{sapp}}
\def\sapp{\par
   \addvspace{4ex}
   \@afterindentfalse
  \secdef\@sapp\@dsapp}
\def\@sapp[#1]#2{\ifnum \c@secnumdepth >\m@ne
        \refstepcounter{sapp}
        \addcontentsline{toc}{sapp}{\thesapp
        \hspace{1em}#1}\else
      \addcontentsline{toc}{sapp}{ #1}\fi
   {\parindent \z@ \raggedright
    \large \bf \thesapp
   \large  \bf \hspace{1em}    #2}\nobreak
   \vskip 4ex   \noindent
\def\theequation{\Alph{app}.\arabic{equation}}}
\def\@dsapp#1{%
{\parindent \z@ \raggedright  \bf #1}\par\nobreak}
\def\l@sapp#1#2{\addpenalty{\@secpenalty}%
   \begingroup
   \@tempdima 3em
     \parindent \z@ \rightskip \@pnumwidth
     \parfillskip -\@pnumwidth
     { \hspace{1em}
     \leavevmode
     #1 \hfil \dotfill \hbox to\@pnumwidth{\hss #2}}\par \nobreak
     \endgroup}

\def\titlepage{\@restonecolfalse\if@twocolumn\@restonecoltrue\onecolumn
     \else \newpage \fi \thispagestyle{empty}\c@page\z@
        \def\thefootnote{\fnsymbol{footnote}} }

\def\endtitlepage{\if@restonecol\twocolumn \else  \fi
        \def\thefootnote{\arabic{footnote}}
        \setcounter{footnote}{0}}  %\c@footnote\z@ }
\relax

\hybrid

\newtoks\@stequation

\def\subequations{\refstepcounter{equation}%
  \edef\@savedequation{\the\c@equation}%
  \@stequation=\expandafter{\theequation}%   %only want \theequation
  \edef\@savedtheequation{\the\@stequation}% %expanded once
  \edef\oldtheequation{\theequation}%
  \setcounter{equation}{0}%
  \def\theequation{\oldtheequation\alph{equation}}}

\def\endsubequations{%
  \setcounter{equation}{\@savedequation}%
  \@stequation=\expandafter{\@savedtheequation}%
  \edef\theequation{\the\@stequation}%
  \global\@ignoretrue}

\makeatletter
\newdimen\normalarrayskip              % skip between lines
\newdimen\minarrayskip                 % minimal skip between lines
\normalarrayskip\baselineskip
\minarrayskip\jot
\newif\ifold             \oldtrue            \def\new{\oldfalse}
\def\arraymode{\ifold\relax\else\displaystyle\fi} % mode of array enrties
\def\eqnumphantom{\phantom{(\theequation)}}     % right phantom in eqnarray
\def\@arrayskip{\ifold\baselineskip\z@\lineskip\z@
     \else
     \baselineskip\minarrayskip\lineskip1\baselineskip\fi}

\def\@arrayclassz{\ifcase \@lastchclass \@acolampacol \or
\@ampacol \or \or \or \@addamp \or
   \@acolampacol \or \@firstampfalse \@acol \fi
\edef\@preamble{\@preamble
  \ifcase \@chnum
     \hfil$\relax\arraymode\@sharp$\hfil
     \or $\relax\arraymode\@sharp$\hfil
     \or \hfil$\relax\arraymode\@sharp$\fi}}

\def\@array[#1]#2{\setbox\@arstrutbox=\hbox{\vrule
     height\arraystretch \ht\strutbox
     depth\arraystretch \dp\strutbox
     width\z@}\@mkpream{#2}\edef\@preamble{\halign \noexpand\@halignto
\bgroup \tabskip\z@ \@arstrut \@preamble \tabskip\z@ \cr}%
\let\@startpbox\@@startpbox \let\@endpbox\@@endpbox
  \if #1t\vtop \else \if#1b\vbox \else \vcenter \fi\fi
  \bgroup \let\par\relax
  \let\@sharp##\let\protect\relax
  \@arrayskip\@preamble}
\def\eqnarray{\stepcounter{equation}%
              \let\@currentlabel=\theequation
              \global\@eqnswtrue
              \global\@eqcnt\z@
              \tabskip\@centering                      %formulae centering
              \let\\=\@eqncr
              $$%
            \halign to \displaywidth  \bgroup
             \eqnumphantom \@eqnsel
      \hskip\@centering                               %right tab%
    $\displaystyle  \tabskip\z@ {##}$%
    &\global\@eqcnt\@ne \hskip 2\arraycolsep
         $ \displaystyle  \arraymode{##}$\hfil
    &\global\@eqcnt\tw@ \hskip 2\arraycolsep
         $\displaystyle\tabskip\z@{##}$\hfil
         \tabskip\@centering
    &{##}\tabskip\z@\cr}
\makeatother

\newtheorem{theor}{Theorem}%[section]      %Usage: {\th Statement}
\newtheorem{prop}{Proposition}[section]           %  ETC ...
\newtheorem{cor}[prop]{Corollary}%[section]
\newtheorem{lem}[prop]{Lemma}%[section]

\def\bea{\begin{eqnarray}}
\def\eea{\end{eqnarray}}

\def\beq{\begin{equation}}
\def\eeq{\end{equation}}
\def\be{\beq\new\begin{array}{c}}
\def\ee{\end{array}\eeq}
\def\bse{\begin{subequations}}                %%%SUBEQUATIONS
\def\ese{\end{subequations}}                 %
\newcommand{\CC}{{\mathbb C}}

\newcommand{\ZZ}{{\mathbb Z}}
\let\z=\ZZ
\def\stackreb#1#2{\mathrel{\mathop{#2}\limits_{#1}}}

\def\al{\alpha}
\def\ve{\varepsilon}
\def\xp{e}
\def\xm{f}
\def\<{\langle}
\def\>{\rangle}
\def\r#1{(\ref{#1})}
\def\ot{\otimes}
\def\a{\alpha}  \def\d{\delta}  
\def\b{\beta} 
\def\D{\Delta}
\newcommand{\Uqg}{{U_{q}^{}(\widehat{\mathfrak{g}})}}
\newcommand{\Uqdva}{{U_{q}^{}(\widehat{\mathfrak{sl}}_{2})}}
\newcommand{\Uqtri}{U_{q}(\widehat{\mathfrak{sl}}_3)}
\let\bn=\be
\let\ed=\ee
\let\rf=\r
\def\sk#1{\left(#1\right)}
\def\Pfpm{{P}^\pm}
\def\Pfp{{P}^+}
\def\Pfm{{P}^-}
\def\res#1{\stackreb{#1}{\mbox{\rm res}}}
\def\si{\sigma}
\def\ant{a}
\def\coun{\varepsilon}
\newcommand{\ein}[2]{e_{#1}[#2]}
\newcommand{\fin}[2]{f_{#1}[#2]}
\newcommand{\hin}[2]{h_{#1}[#2]}
\def\U{\overline{U}}
\def\i{\iota}
\def\l{\lambda}
\newcommand{\UqtriD}{{\U_{q}(\widehat{\mathfrak{sl}}_{3})}}
\def\op{{{\hat{\mathbf \omega}}}}

\begin{document}
\thispagestyle{empty}
\setcounter{page}1

\begin{center}
\hfill ITEP-TH-56/04\\
\hfill MPIM 2005-16\\
\hfill math.QA/0610433\\
\bigskip\bigskip
{\Large\bf Weight function for the quantum\\ affine algebra  $\Uqtri$}\\
\bigskip
\bigskip
{\bf
Sergey Khoroshkin$^{\circ,\diamond}$\footnote{E-mail: khor@itep.ru},
Stanislav Pakuliak$^{\circ,\star,\diamond}$\footnote{E-mail: pakuliak@theor.jinr.ru}
}\\
\bigskip
$^\circ${\it Institute of Theoretical \ and Experimental Physics, 117259
Moscow, Russia}\\
$^\star${\it Bogoliubov Laboratory of Theoretical Physics, JINR,
141980 Dubna, Moscow region, Russia}\\
$^\diamond${\it Max-Plank-Institut f\"ur Mathematik, Vivatsgasse 7,
D-53111 Bonn, Germany}
\\
\bigskip
\bigskip
%{Revised \today}
\end{center}
\begin{abstract} We give a precise expression for the universal weight function
 of the quantum affine algebra $\Uqtri$.
The calculations use the technique of projecting products of Drinfeld currents
on the intersections of Borel subalgebras.
\end{abstract}
\footnotesize
%\tableofcontents
%\newpage
\normalsize

\section{Introduction}
 The ideology of a nested Bethe ansatz \cite{KR83} prescribes two steps
 for describing the transfer-matrix eigenvectors in finite-dimensional
 representations of a quantum affine algebra. First, specific rational
 functions with values in the representation should be constructed;
 second, a system of Bethe equations for these functions should be solved.
 These rational vector-valued functions are called off-shell (nested)
 Bethe vectors. They can serve as a generating system of vectors of a
 finite-dimensional representation of a quantum affine algebra. We use
 the equivalent name 'weight function', which came from applications in
 difference Knizhnik-Zamolodchikov equations \cite{S}, \cite{VT}.

 A general construction of a weight function for a quantum affine algebra
 $\Uqg$ was recently suggested \cite{EKPT}. This construction uses the
 existence of two different types of Borel subalgebras in a quantum affine
 algebra. One type is related to the realization of $\Uqg$ as a quantized
 Kac-Moody algebra, and the other comes from the current realization of $\Uqg$
 proposed by Drinfeld \cite{D88}. The weight function is defined as the
 projection of a product of Drinfeld currents on the intersection of Borel
 subalgebras of $\Uqg$ of different
 types (see Sec.~\ref{universal-w-f}).

 Our goal in this paper is to develop a technique for calculation the weight
 function starting from the definition in \cite{EKPT}.
 According to this definition, for calculation the weight function, the product
 of Drinfeld currents must be arranged in a normally ordered form.
 Then only those terms are kept that belong to the intersection of Borel
 subalgebras of different types. The normal-ordering procedure requires
 investigating the current adjoint action and the composed root currents,
 introduced in \cite{DKh}. The final result is a precise universal
 expression for the weight function of $\Uqtri$, which can then be specialized
 to any finite-dimensional representation of $\Uqtri$.
 %We note that the calculation
 %of the weight function for $U_q({\widehat{\mathfrak gl}}_N)$ at the level
 %of a tensor product of evaluation representations can be found in \cite{T}.

 This paper is organized as follows.
 In Section~\ref{sl3-intro}, we introduce the main objects of the investigation.
 Section~\ref{main-results} is devoted to formulating the main results. They contain
 a precise expression for the weight function of $\Uqtri$ (Theorems~\ref{act-cal}
 and \ref{prstring}). As a particular case, we give an expression for the weight
 function of $\Uqdva$ in an integral form (Theorem~\ref{prstring}). The kernel of
 the integral is a well-known partition function, which coincides with the partition
 function of the six-vertex model on a finite square lattice with fixed domain-wall
 boundary conditions. Later, we need a combinatorial identity for this kernel, which
 we prove by observing the self-adjointness of the projection operators
 (Proposition~\ref{PPprop}).

 Sections~\ref{anal-st} and \ref{adjoint-act} are devoted to proving
 the main statements, which includes studying the analytic properties of
 composed currents and related products of currents (strings) and of current
 adjoint actions. We also note an important role of symmetrization procedures,
 based on the properties of the analytic continuation of the products of
 currents and of their projections (see Proposition~\ref{anprop}).
 In the appendices, we give the necessary properties of the opposite
 projection operator, commutation relations between currents and their
 projections, and another proof of the main result.

\section{Basic notation}\label{sl3-intro}
\subsection{ $\Uqtri$ in Chevalley generators}\label{ch-gen}

The quantum affine algebra $\Uqtri$
 is generated by Chevalley\footnote{In what follows we do not
use a grading operator and set the central charge equal to zero.
Such an algebra is usually denoted by $U'_q(\widetilde{{\mathfrak
sl}}_3)$.} generators $e_{\pm\a_i}$ and $k^{\pm1}_{\a_i}$, where
$i=0,1,2$ and $\prod_{i=0}^2 k_{\a_i}=1$, subject to the  relations
\be\label{Chevalley}
%[d,e_{\pm\al_i}]\,=\,\pm\,\delta_{i0}e_{\pm\al_i},\quad
k_{\al_i}e_{\pm\al_j}k^{-1}_{\al_i}\,=\,q_i^{\pm a_{ij}}e_{\pm\al_j}\,,
\quad
[e_{\al_i},e_{-\al_j}]\,=\,
\delta_{ij}\frac{k_{\al_i}-k^{-1}_{\al_i}}{q_i-q_i^{-1}}\,,
\ee
\be\label{re2}
e^2_{\pm\al_i}e_{\pm\al_j}+[2]_qe_{\pm\al_i}e_{\pm\al_j}e_{\pm\al_i}+
e_{\pm\al_j}e^2_{\pm\al_i}
=0\,,\quad
i\neq j\,,\quad (\a_i,\a_j)=-1\,,
\ee
where $[n]_q=\frac{q^n-q^{-n}}{q-q^{-1}}$ is the Gauss $q$-number and
$a_{ij}=(\al_i,\al_j)$ is symmetrized
Cartan matrix of the affine algebra $\widehat{\mathfrak{sl}}_3$,
\begin{equation}\new\begin{array}{c}\label{Cartan}
a_{ij}=(\al_i,\al_j)=\sk{\begin{array}{ccc}2&-1&-1\\
-1&2&-1\\-1&-1&2\end{array}}.
\end{array}\end{equation}

One of the possible
Hopf structures (which we call the standard Hopf structure)
is given by the formulas
\be\label{copr}
\Delta(e_{\al_i})=e_{\al_i}\otimes 1+k_{\al_i}\otimes e_{\al_i}\,,\quad
\Delta(e_{-\al_i})=1\otimes e_{-\al_i}+e_{-\al_i}\otimes
k^{-1}_{\al_i}\,,\\
\Delta(k_{\al_i})\,=\,k_{\al_i}\otimes k_{\al_i}\,,\quad
\coun(e_{\pm\al_i})=0\,,\quad  \coun(k^{\pm1}_{\al_i})=1\,,\\
\ant(e_{\al_i})=-k^{-1}_{\al_i}e_{\al_i}\,,\quad
\ant(e_{-\al_i})=-e_{-\al_i}k_{\al_i}\,,\quad
\ant(k^{\pm1}_{\al_i}) =  k^{\mp1}_{\al_i}\,,
\ee
where $\Delta$, $\coun$ and $\ant$ are the respective comultiplication, counit,
and antipode maps.

\subsection{Current realization of the algebra
$\Uqtri$}
As does any quantum affine algebra, $\Uqtri$ admits a current realization
\cite{D88}.  In this description
 (we again assume that the central charge is zero),
 $\Uqtri$ is generated by the elements
$\ein{i}{n}$ and $\fin{i}{n}$, where $i=\a,\b$ and $n\in\ZZ$, and $k_i^{\pm 1}$
and $\hin{i}{n}$, where  $i=\a,\b$ and $n\in\ZZ\,%\setminus
\backslash \{0\}$.
They are gathered in the generating functions
\be\label{currents}
e_i(z)\,=\,\sum_{n\in\z}\ein{i}{n}z^{-n}\ ,
\quad f_i(z)\,=\,\sum_{n\in\z}\fin{i}{n}z^{-n}\ , \\
\psi^\pm_i(z)\,=\sum_{n>0}\psi_i^{\pm}[n]z^{\mp n}
=k_{i}^{\pm1}
\,
\exp\sk{\pm(q-q^{-1})\sum_{n>0}\hin{i}{\pm n}z^{\mp n}}\ ,
\ee
which satisfy the  relations
\bn
\begin{array}{rcl}
(z-q^{({i},{j})}w)e_{i}(z)e_{j}(w)&=&
 e_{j}(w)e_{i}(z)(q^{({i},{j})}z-w)\ ,\\
(z-q^{-({i},{j})}w)f_{i}(z)f_{j}(w)&=&
 f_{j}(w)f_{i}(z)(q^{-({i},{j})}z-w)\ ,
\end{array}
\label{1}
\ed
\bn
\begin{array}{rcl}
\ds
 \psi_{i}^\pm(z)e_{j}(w)\left(\psi_{i}^\pm(z)\right)^{-1}&=&\ds
\frac{(q^{({i},{j})}z-w)}
{(z- q^{({i},{j})}w)}e_{j}(w)\ ,\\
\ds
\psi_{i}^\pm(z)f_{j}(w)\left(\psi_{i}^\pm(z)\right)^{-1}&=&\ds
\frac{(q^{-({i},{j})}z-w)}
{(z- q^{-({i},{j})}w)}f_{j}(w)\ ,
\end{array}
\label{3}
\ed
\bn
\psi_{i}^\mu(z)\psi_{j}^\nu(w)=
 \psi_{j}^\nu(w)\psi_{i}^\mu(z)\ ,\quad \mu,\nu=\pm\,,
\label{7a}
\ed
\bn
[e_{i}(z),f_{j}(w)]=
 \frac{\delta_{ij}\delta(z/w)}{q_{}-q^{-1}_{}}\left(
\psi^+_{i}(z)-\psi^-_{i}(w)\right)\,,
\label{10}
\ed
where $i,j=\a,\b$, $\delta(z)=\sum_{k\in \z}z^k$,
$(\a,\a)=(\b,\b)=2$, $(\a,\b)=-1$,
and
%\bn
\begin{eqnarray}\label{serre2}%{rcl}
\ds \mathop{\rm Sym}\limits_{z_1,z_2}
\sk{e_i(z_1)e_i(z_2)e_j(w) -(q+q^{-1})e_i(z_1)e_j(w)e_i(z_2)
+e_j(w)e_i(z_1)e_i(z_2)  }&=&0\,,\\ \label{serre1}
\ds \mathop{\rm Sym}\limits_{z_1,z_2}
\sk{f_i(z_1)f_i(z_2)f_j(w) -(q+q^{-1})f_i(z_1)f_j(w)f_i(z_2)
+f_j(w)f_i(z_1)f_i(z_2)  }&=&0\,,
\end{eqnarray}
%\ed
where  $i,j=\a,\b$, $i\not=j$. The assignment
\bn
\begin{array}{ll}
\ds k_{\a_1}^{}\mapsto k_\a^{}\,,\qquad
k_{\a_2}^{}\mapsto k_\b^{}\,,\qquad
&k_{\a_0}^{}\mapsto k_\a^{- 1}k_\b^{- 1},\\
\ds e_{\a_1}\mapsto \ein{\a}{0}\,,\qquad
e_{\a_2}\mapsto \ein{\b}{0}\,,\qquad
&e_{-\a_1}\mapsto \fin{\a}{0}\,,\qquad
e_{-\a_2}\mapsto \fin{\b}{0}\,,\\
\ds e_{\a_0}\mapsto \fin{\b}{1}\fin{\a}{0}-q\fin{\a}{0}\fin{\b}{1}\,,\qquad
&e_{-\a_0}\mapsto \ein{\a}{0}\ein{\b}{-1}-q^{-1}\ein{\b}{-1}\ein{\a}{0}
\end{array}
\ed
establishes the isomorphism of the two realizations.

The algebra $\Uqtri$ admits a natural completion
 $\UqtriD= U^{(D)}_q(\widehat{\mathfrak{sl}}_{3})$, which
can be described as the minimal extension of  $\Uqtri$ and which
acts in all representations of  $\Uqtri$ that are highest-weight
representations of $U_q(\mathfrak{b}_+)$
% (in the category $\mathcal{O}$ of the highest weight representations).
 (see Sec. 2.2 in \cite{DKhP1} for the details).

In a highest-weight representation of $\Uqtri$, any matrix
coefficients
 of an arbitrary product of the currents $a_1(z_1)...,a_n(z_n)$ are formal power
series in the space
$$
 {\bf C}[z_1,z_1^{-1},...,z_m,z_m^{-1}]\left[\left[\frac{z_2}{z_1},
\frac{z_3}{z_2},...,\frac{z_m}{z_{m-1}}\right]\right]
$$
and converge to a rational function in the domain
$|z_1|\gg|z_2|\gg\cdots\gg|z_m|$ (see \cite{E}, \cite{DKhP}). This
observation and commutation relations \rf{1}, which dictate the rule
for the analytic continuation from the above domain, allow
considering products of currents as meromorphic functions with
values in $\UqtriD$. In what follows, we freely use this analytic
language and replace formal integrals with contour integrals in this
formalism. An integral without the contour specified always means a
formal integral.

Another Hopf structure in $\Uqtri$ is naturally related to the
current realization. In terms of currents, it is given by
\bse\label{du1} \be\label{du1a} \Delta^{(D)}\xp_i(z)=\xp_i(z)\otimes
1+\psi_i^-(z)\otimes \xp_i(z)\, , \ee \be\label{du1b}
\Delta^{(D)}\xm_i(z)=1\otimes \xm_i(z)+\xm_i(z)\otimes \psi_i^+(z)\,
, \ee \be\label{du1c} \Delta^{(D)}\psi_i^\pm(z)=\psi_i^+(z)\otimes
\psi_i^\pm(z)\, , \ee \be\label{du1d}
\ant\sk{e_i(z)}=-\sk{\psi_i^-\sk{z}}^{-1}\, e_i(z)\,,\quad
\ant\sk{f_i(z)}=-f_i(z)\, \sk{\psi_i^+\sk{z}}^{-1}, \ee
\be\label{du1f} \ant\sk{\psi_i^\pm(z)}=\sk{\psi_i^\pm(z)}^{-1},\quad
\coun(e_i(z))=\coun(f_i(z))=0\,,\quad \coun(\psi_i^\pm(z))=1\, . \ee
\ese The comultiplications $\Delta$ in Sec.~\ref{ch-gen} and
$\Delta^{(D)}$ are related by the twist, which can be described
explicitly (see \cite{DKhP1}).

\subsection{Borel subalgebras of $\Uqtri$}\label{bor-sub}

We let $U_q(\mathfrak{b}_+)$ denote the subalgebra of $\Uqtri$ generated by
the elements $e_{\a_i}$ and $k_{\a_i}^{\pm 1}$, $i=0,1,2$. We also let
$U_q(\mathfrak{b}_-)$ denote the subalgebra of $\Uqtri$ generated by
the elements $e_{-\a_i}$ and $k_{\a_i}^{\pm 1}$, $i=0,1,2$.

The algebras  $U_q(\mathfrak{b}_\pm)$ are Hopf subalgebras of $\Uqtri$
with respect to the standard comultiplication $\Delta$ and serve as
$q$-deformations of the enveloping algebras of opposite Borel subalgebras of
the Lie algebra $\widehat{\mathfrak{sl}}_{3}$. We call them the standard Borel
subalgebras. They contain subalgebras $U_q(\mathfrak{n}_\pm)$ generated by the
elements $e_{\pm \a_i}$, $i=0,1,2$.

 The subalgebra
$U_q(\mathfrak{n}_+)$ is a left coideal of $U_q(\mathfrak{b}_+)$ with respect
to the standard comultiplication, and the subalgebra $U_q(\mathfrak{n}_-)$ is
a  right coideal of $U_q(\mathfrak{b}_-)$ with respect
to the standard comultiplication, i.e.,
$$\Delta(U_q(\mathfrak{n}_+))\subset U_q(\mathfrak{b}_+)\otimes
U_q(\mathfrak{n}_+)\,, \qquad
\Delta(U_q(\mathfrak{n}_-))\subset U_q(\mathfrak{n}_-))\otimes
U_q(\mathfrak{b}_-) \, .
$$

Borel subalgebras of another type are related to the current realization of
$\Uqtri$. We let $U_F$ denote the subalgebra of $\Uqtri$ generated by the elements
$k_i^{\pm 1}$ and $\fin{i}{n}$, where $i=\a,\b$ and $n\in\ZZ$, and $\hin{i}{n}$,
where $i=\a,\b$ and  $n>0$. Its completion $\U_F$ is a Hopf subalgebra of $\UqtriD$
with respect to the comultiplication $\Delta^{(D)}$. We call $U_F$ the current Borel
subalgebra. It contains the subalgebra $U_f$ generated by the elements
$\fin{i}{n}$, where $i=\a,\b$ and  $n\in\ZZ$.
The completed algebra $\U_f$ is a right coideal of $\U_F$ with respect to the
comultiplication $\Delta^{(D)}$ and serves as a $q$-deformed enveloping algebra
 of the algebra of currents valued in $\mathfrak{n}_-$.

 The opposite current  Borel subalgebra $U_E$ is
 generated by the elements
$k_i^{\pm 1}$ and $\ein{i}{n}$, where $i=\a,\b$, $n\in\ZZ$  and
 by the elements $\hin{i}{n}$,
$i=\a,\b$, $n<0$.

\subsection{Projections $P^\pm$ on intersections of Borel subalgebras}

We let $U_F^+$ and $U_f^-$ denote the subalgebras of the current
Borel algebra $U_F$, \be\label{F-dec} U_f^-=U_F\cap
U_q(\mathfrak{n}_-)\,,\qquad U_F^+=U_F\cap U_q(\mathfrak{b}_+)\,.
\ee {}For any $x\in \Uqtri$, we let ${\rm ad}_x: \Uqtri\to \Uqtri$
be the operator of the adjoint action of $x$ in $\Uqtri$. It is
defined by the relation
$${\rm ad}_x(y)=\sum_ja(x'_j)\cdot y\cdot x''_j, \qquad {\rm if}\qquad
\Delta(x)=\sum_jx'_j\otimes x''_j\ .$$
{}For $i=\a,\b$, let $S_i$ be the operator ${\rm ad}_{\fin{i}{0}}$ such that
\bn\label{coad}
S_{i}(y)=y\fin{i}{0}-\fin{i}{0}k_i^{}yk_i^{-1}.
\ed
We call $S_i$ the {\it screening operators}.
\begin{prop}\label{subalgebras}${}$
\begin{itemize}
\item[(i)]
 The algebra  $U_f^-$ is generated by the elements
$\fin{i}{n}$, where $i=\a,\b$ and $ n\leq0$; the algebra $U_F^+$
is generated by the elements $k_i^{\pm 1}$, $\fin{i}{n}$, and $\hin{i}{n}$,
 where $i=\a,\b$ and $ n>0$, and by the element
\bn\label{fab}
\fin{\a+\b}{1}=\fin{\b}{1}\fin{\a}{0}-q\fin{\a}{0}\fin{\b}{1}=-\left(
\fin{\a}{1}\fin{\b}{0}-q\fin{\b}{0}\fin{\a}{1}\right)\, .\ed
\item[(ii)] The subalgebras $U_f^-$ and  $U_F^+$ are
  invariant under the action of the screening operators
   $S_{i}$, $i=\a,\b$.
\item[(iii)] The subalgebra $U_F^+$ is a left coideal of $U_F^{}$
with respect to the comultiplication $\Delta^{(D)}$;
the subalgebra  $U_f^-$ is a right coideal of $U_F^{}$
with respect to the comultiplication $\Delta^{(D)}$.
\item[(iv)] The multiplication in $U_F$ establishes an isomorphism of the vector
spaces $U_F$ and $U_f^-\ot U_F^+$\,.
\end{itemize}
\end{prop}
\noindent
{\it Proof}.\  Statement (ii) can be verified as follows:
$$
\begin{array}{l}
\ds
S_{\a}\sk{\fin{\b}{1}\fin{\a}{0}-q\fin{\a}{0}\fin{\b}{1}}
=\\[3mm]
\ds\qquad = \sk{\fin{\b}{1}\fin{\a}{0}-q\fin{\a}{0}\fin{\b}{1}}\fin{\a}{0}-q^{-1}
\fin{\a}{0}
\sk{\fin{\b}{1}\fin{\a}{0}-q\fin{\a}{0}\fin{\b}{1}}=\\[3mm]
\ds\qquad =
\fin{\b}{1}\fin{\a}{0}\fin{\a}{0}-(q+q^{-1})\fin{\a}{0}\fin{\b}{1}
\fin{\a}{0}+\fin{\a}{0}\fin{\a}{0}\fin{\b}{1}=0\,.
\end{array}
$$
To prove (iii), we use formula \r{du1a} written in terms of modes as
$$
\Delta^{(D)}\sk{f_i[n]}=1\otimes f_i[n]+\sum_{k\geq0}
f_i[n-k]\otimes \psi^+_i[k]\, .
$$
We must show that $\Delta^{(D)}\sk{\fin{\b}{1}\fin{\a}{0}-q\fin{\a}{0}\fin{\b}{1}}
\in U_F\otimes U_F^+$. The formula above shows that it suffices to verify
that $\psi^+_\b[k]\fin{\a}{0}-q\fin{\a}{0}\psi^+_\b[k]  \in U_F^+$. But this
holds because of the relation
$$
\psi^+_\b(z)\fin{\a}{0}-q\fin{\a}{0}\psi^+_\b(z)=
(q^{2}-1)
\sum_{n=1}^\infty
(qz)^{-n}\psi^+_{\b}(z)\fin{\a}{n}\,.
$$
\hfill$\Box$

We define the operators $P=P^+$: $U_F\to U_F^+$ and $P^-$:  $U_F\to
U_f^-$ by the relations \bn \label{Pdef} P(f_1f_2)\ =\
P^+(f_1f_2)=\varepsilon(f_1)f_2, \qquad
P^-(f_1f_2)=f_1\varepsilon(f_2)\ \ed for any $f_1\in U_f^-$ and
$f_2\in U_F^+$. Proposition \ref{subalgebras} implies that the
algebras $U_f^-$ and $U_F^+$ satisfy conditions {\it (i)} and  {\it
(ii)} in Sec. 4.1 in \cite{DKhP} (also see Section 6 in \cite{ER})
with respect to the comultiplication
$\left(\Delta^{(D)}\right)^{op}$. By \cite{DKhP}, the operators
$P^\pm$ are therefore well-defined projection operators,
$\left(P^\pm\right)^2= P^\pm$, which admit extensions to the
completed algebra $\U_F$ such that for any $f\in \U_F$, the
canonical decomposition \be\label{addF} f=\sum_iP^-(f''_i)\cdot
P^+(f'_i)\,,\quad {\rm if}\quad \Delta^{(D)}(f)=\sum_if'_i\otimes
f''_i\, \ee holds. We call expressions of the form
$f=\sum_i\tilde{f}_i{\tilde{f}}'_i$, where $\tilde{f}_i\in \U_f^-$
and ${\tilde{f}}'_i\in \U_F^+$, the {\it normally ordered
expansion}. The normally ordered expansion is compatible with the
action of the algebra
 $\Uqtri$ in highest-weight representations.
% представлениях категории $\mathcal{O}$.
Expression
\rf{addF} gives an ordered expansion of an arbitrary element
 $f\in \overline{U}_F$.

\subsection{Composite current and strings}\label{com-cs}

We define the generating function of the elements in $\UqtriD$
\be\label{rec-f} f_{\a+\b}(z)= \oint f_{\a}(z) f_{\b}(w)\
\frac{dw}{w}- \oint \frac{q^{-1}-z/w}{1-q^{-1}z/w} f_{\b}(w)
f_{\a}(z)\ \frac{dw}{w}\,, \ed where the formal integral $\oint
g(w)\frac{dw}{w}$ of a Laurent series $g(w)=\sum_{k\in\ZZ}g_kw^{-k}$
means taking its  coefficient $g_{0}$.

We can also write the formal integral in the analytic language \cite{DKh},
\bn \label{rec+}
f_{\a+\b}(z)=-\res{w=zq^{-1}}f_{\a}(z)f_{\b}(w)\,\frac{dw}{w}
\ed
such that the relation
\be\label{com-f}
f_{\a}(z)f_{\b}(w)=
\ds \frac{1-qz/w}{q-z/w}f_{\b}(w)f_{\a}(z)+
\delta(zq^{-1}/w)f_{\a+\b}(z)
\ee
{}holds in the algebra $\UqtriD$. For any $a,b\in\ZZ_{\geq 0}$, the products
\bn\label{string}
f_\a(u_1)\cdots \!f_\a(u_a)f_{\a+\b}(u_{a+1})\cdots \!f_{\a+\b}(u_{a+b}),
\quad
f_{\a+\b}(u_{1})\cdots \!f_{\a+\b}(u_{a})f_\b(u_{a+1})\cdots \!f_\b(u_{a+b})
\ed
are called  {\it strings}.
The products
\bn\label{stringop}
f_{\a+\b}(u_{a+b})\cdots \!f_{\a+\b}(u_{a+1})f_\a(u_a)\cdots \!f_\a(u_1),
\quad
f_\b(u_{a+b})\cdots \!f_\b(u_{a+1})
f_{\a+\b}(u_{a})\cdots \!f_{\a+\b}(u_{1})
\ed
are called {\it opposite strings} to strings \rf{string}.
 The strings have nice analytic properties, which are crucial
for their use in this paper. These properties are listed in
Proposition~\ref{hyper}.

\setcounter{equation}{0}
\section{Main results}\label{main-results}
\subsection{Universal weight function}\label{universal-w-f}
 Let $V$ be a representation of $\Uqtri$ and $v$ be a vector
in $V$. We call $v$ a highest-weight vector with respect to the current Borel
subalgebra $U_E$, if
\bn
\begin{array}{rcl}
e_i(z)v&=&0\,, \\
\psi_i^{\pm}(z)v&=&\l_i(z)v\,, \qquad i=\a,\b\,,
\end{array}
\ed
where $\l_i(z)$ is a meromorphic function decomposed into a series in
$z^{-1}$ for $\psi^+_i(z)$ and into a series in $z$ for $\psi^-_i(z)$.
The representation $V$ is called a representation with the highest-weight vector
$v\in V$ with respect to $U_E$ if it is generated by $v$ over
$\Uqtri$.

Let $\Pi$ denote the two-element set $\{\a,\b\}$ of positive simple
roots of the Lie algebra ${\mathfrak sl}_3$.
An ordered set $I=\{a_1,...,a_{|I|}\}$, together with a map
$\i: I\to \Pi$, is called an ordered $\Pi$-multiset.

We suppose that for any ordered $\Pi$-multiset $I$, $|I|=n$, a formal
series  $W(t_{i_1},...,t_{i_n})$ $\in$ $ U\{t_{i_1},...,t_{i_n}\}$,
$i_k\in I$, is chosen, where
\bn%% W(t_{i_1},...,t_{i_n})\in
 U\{t_{i_1},...,t_{i_n}\}\ = \  \Uqtri \label{Wser}
[t_{i_1},t_{i_1}^{-1},...,t_{i_n},t_{i_n}^{-1}]\left[\left[\frac{t_{i_2}}{t_{i_1}},
 \frac{t_{i_3}}{t_{i_2}},
...,\frac{t_{i_n}}{t_{i_{n-1}}},
\frac{1}{t_{i_n}}\right]\right] ,
\ee

i.e., $W(t_{i_1},...,t_{i_n})$ is a formal power
series in the variables
 ${t_{i_2}}/{t_{i_1}}, {t_{i_3}}/{t_{i_2}}
,...,{t_{i_n}}/{t_{i_{n-1}}},
1/t_{i_n}$ with coefficients in polynomials
$\Uqtri[t_{i_1},t_{i_1}^{-1},...,t_{i_n},t_{i_n}^{-1}]$ such that
\begin{itemize}
\item[1)] for any representation $V$ that is highest-weight with respect to $U_E$
 with the highest weight vector $v$, the function
$$w_V(t_{i_1},...,t_{i_n})=W(t_{i_1},...,t_{i_n})v$$
converges in the domain  $|t_{i_1}|\gg\cdots\gg|t_{i_n}|$ to a meromorphic $V$-valued function,
\item[2)] if $I=\emptyset$, then $W=1$ and $w_V=v$, and
\item[3)] if $V=V_1\otimes V_2$ is a tensor product of highest-weight
representations with the highest-weight vectors $v_1$ and $v_2$ and the
highest-weight series $\{\l_i^{(1)}(z)\}$ and
$\{\l_i^{(2)}(z)\}$, $i=\a,\b$, then for any ordered $\Pi$-multiset $I$, we
have
\bn
\begin{array}{c}
w_{V}(\{t_{a}|_{ a\in\ I}\})=
\sum\limits_{I=I_1\coprod I_2}
w_{V_1}(\{t_{a}|_{ a\in\ I_1}\})
\otimes
w_{V_2}(\{t_{a}|_{ a\in\ I_2}\})
\times \\
\times \prod\limits_{{a\in I_1}}
\l^{(2)}_{\i(a)}(t_{a})
\times
\prod\limits_{{a<b,\ a\in I_1,\ b\in I_2}}
\frac{q^{-(\i(a),\i(b))}t_{a}-t_{b}}
{t_{a}-q^{-(\i(a),\i(b))}t_{b}}\,.
\label{weight1}
\end{array}
\ed
\end{itemize}
A collection $W(t_{i_1},...,t_{i_n})$ is called a {\it
universal weight function}. A collection $w(t_{i_1},...,t_{i_n})$
is called a {\it weight function}.

If solutions of the corresponding system of Bethe equations are
chosen as parameters in the weight function, then a set of Bethe
vectors is obtained. The weight function with free parameters was
systematically used to investigate solutions of $q$-difference
Knizhnick-Zamolodchikov equations %\cite{EKPT},
\cite{S}, \cite{VT}.

Let $I=\{{i_1},...,{i_n}\}$ be an ordered $\Pi$-multiset. We set
\bn \label{Wt}
W(t_{i_1},...,t_{i_n})=P\left(f_{\i(i_1)}(t_{i_1})\cdots f_{\i(i_n)}(t_{i_n})
\right).
\ed
The main result in paper \cite{EKPT} can be formulated as follows
in the particular case of $\Uqtri$.
\medskip

%%\begin{theor}%\nonumber
{\bf Theorem.\, }%%\label{weightf}
 \cite{EKPT} {\it The collection $W(t_{i_1},...,t_{i_n})$ defined in
\rf{Wt} is a universal weight function .}
\medskip

We again note that all the expressions for universal weight
functions $W(t_{i_1},...,t_{i_n})$ are to be understood as formal
series in the variables ${t_{i_2}}/{t_{i_1}}, {t_{i_3}}/{t_{i_2}}
,...,{t_{i_n}}/{t_{i_{n-1}}}, 1/t_{i_n}$. If we deal with a weight
function $w(t_{i_1},...,t_{i_n})$ that is a vector-valued rational
function, then there is no difference in the choice of the domain
where this function is expanded (see Sec.~\ref{razdel5.1} for more
details).

\subsection{Reduction to projections of strings}

Let $S_n$ be the group of permutations on $n$ elements. For any  set
$\overline{t}=\{t_{1},...,t_{n}\}$ of variables $t_{1},...,t_{n}$
 and any $\sigma
\in $ $S_n$, we let $^\sigma\overline{t}$ denote the set
$\{t_{\sigma(1)},..., t_{\sigma(n)}\}$. We keep the symbol
$\hat{\omega}$ for the longest element of the group $S_n$. In this
notation, the set $^\op\overline{t}$ means the set $\overline{t}$
with the reversed order: $^\op\overline{t}$ $=$
$\{t_{n},...,t_{1}\}$.

The group $S_n$ acts naturally
in the space of vector-valued meromorphic functions of $n$ variables
$\overline{t}=\{t_1,...,t_n\}$
by the rule $F(\overline{t})\mapsto {}^{\sigma^{}} F(\overline{t})$, where
$$^\sigma F(\overline{t})= {}^\sigma F(t_1,...,t_n)=
F(t_{\sigma(1)},...,t_{\sigma(n)}) =F(^\sigma\overline{t})\, .$$
We now suppose that  $F(\overline{t})$
 is a series  in the domain $|t_1|\gg\cdots\gg|t_n|$ with values in the vector space $V$,
 i.e.,  $F(\overline{t})$ belongs to the space
\be \label{series}
V[t_{1},t_{1}^{-1},...,t_{n},t_{n}^{-1}]\left[\left[\frac{t_{2}}{t_{1}},
 \frac{t_{3}}{t_{2}},
...,\frac{t_{n}}{t_{{n-1}}}, \frac{1}{t_{n}}\right]\right] . \ee We
suppose that this series converges in the domain
$|t_1|\gg\cdots\gg|t_n|$ to an analytic function and that for any
$\sigma\in S_n$, this analytic function admits an analytic
continuation to the domain
$|t_{\sigma(1)}|\gg\cdots\gg|t_{\sigma(n)}|$. We then set ${}^\sigma
F(\overline{t})$ equal to the formal  series representing the
analytic continuation of the function
 $F(t_{\sigma(1)},...,t_{\sigma(n)})$
to the domain
$|t_1|\gg\cdots\gg|t_n|$. Therefore,  ${}^\sigma F(\overline{t})$ is again
a series in (3.5).

With this convention, the symmetrization
${\rm Sym}_{t}^{n}\ F(t_1,...,t_n)$ of a function $F(t_1,...,t_n)$,
as well as of a series $F(t_1,...,t_n)$ in a domain $|t_1|\gg\cdots\gg|t_n|$, is
the sum ${\rm Sym}_{t}^{n}\ F(t_1,...,t_n)=\sum_{\sigma\in S_n}{}^\sigma
F(t_1,...,t_n)$.
The $q$-symmetrization of a function $F(\overline{t})$ of $n$ variables
or of a  series $F(\overline{t})$ in a domain $|t_1|\gg\cdots\gg|t_n|$ is
defined as
\be\label{qs1ra}
\overline{\rm Sym}_{t}^{n}\ F(\overline{t})=
\sum_{\si\in S_{n}}
\prod_{\ell<\ell'\atop \si^{}(\ell)>\si^{}(\ell')}
\frac{q^{-1}-qt_{\si(\ell)}/t_{\si(\ell')}}
{q-q^{-1}t_{\si(\ell)}/t_{\si(\ell')}}\
{}^\si F(\overline{ t})\,.
\ee
Symmetrization of a series that is convergent in a different asymptotic zone
is defined in analogously.

Universal weight function \rf{Wt} allows analytic continuations
to different asymptotic zones because the operator $P$ extends to a projection
operator in the completed algebra $\overline{U}_F$,
where the analytic continuation of the products of currents is well defined.

Let $I=\{i_1,...,i_n\}$ be an ordered $\Pi$-multiset. For any
permutation $\sigma\in S_n$, we let $^\sigma I$ denote an ordered
 $\Pi$-multiset $^\sigma I=\{i_{\sigma(1)},...,i_{\sigma(n)}\}$
 that differs from $I$ by permutations of the elements but has the same
 map $\iota: I\to \Pi$. Let $W(t_{i_{\sigma(1)}}\cdots t_{i_{\sigma(n)}})$ be
 a universal weight function corresponding to the ordered set
 $^\sigma I$ and $\widetilde{W}\sk{t_{i_1}\cdots t_{i_n}}$ be the analytic
 continuation of the weight function ${W}\sk{t_{i_1}\cdots t_{i_n}}$ to the domain
 $|t_{i_{\sigma(1)}}|\gg ... \gg |t_{i_{\sigma(n)}}|$.

\begin{prop}\label{symmetry}
Universal weight function  \rf{Wt} satisfies the relations
\begin{equation}\label{symmetry2}
 W(t_{i_{\sigma(1)}}\cdots t_{i_{\sigma(n)}})=
\prod\limits_{k<l\atop \si^{-1}(k)>\si^{-1}(l)}
\frac{q^{(\i(i_k),\ \i(i_l))} -
\frac{\ds t_{i_l}}{\ds t_{i_k}}}
{1-q^{(\i(i_k),\ \i(i_l))}
\frac{\ds t_{i_l}}{\ds t_{i_k}}}
 \widetilde{W}\sk{t_1\cdots t_n}\,.
\end{equation}
\end{prop}
This proposition is a direct consequence of
Proposition~\ref{anprop}. It follows from Proposition~\ref{symmetry}
that the universal weight function for $\Uqtri$ is completely
defined  by the expression \bn\label{Wts} W(t_1,...,
t_a,s_1,...,s_b)\ =\ P(f_\a(t_1)\cdots f_\a(t_a)f_\b(s_1)\cdots
f_\b(s_b))\, . \ed In this paper, we suggest an explicit expression
for  function \rf{Wts} in terms of the current generators of
$\Uqtri$.

For the sets of variables $\bar t=\{t_1,...,t_k\}$ and $\bar s=\{s_1,...,s_k\}$,
we define the series
\be\label{rat-Y}
\begin{array}{rcl}
\ds Y(\bar t;\bar s)&=&\ds \prod_{i=1}^k\frac{1}{1-s_i/t_i}
\prod_{j=1}^{i-1}\frac{q-q^{-1}s_j/t_i}{1-s_j/t_i}\
%%\\ &=&\ds
\ = \ \prod_{i=1}^k\frac{1}{1-s_i/t_i}
\prod_{j=i+1}^{k}\frac{q-q^{-1}s_i/t_j}{1-s_i/t_j} \,, \\
\ds Z(\bar t;\bar s)&=&Y(\bar t;\bar s) \prod_{i=1}^k
\frac{s_i}{t_i}
\end {array}
\ee

\begin{theor}\label{act-cal}
Universal weight function \rf{Wts} can be written as
\be\label{po1}
\begin{array}{l}
\ds W(t_1,..., t_a,s_1,...,s_b)\ =\
P\sk{f_\a(t_1)\cdots f_\a(t_{a})f_\b(s_1)\cdots
f_\b(s_{b})}=\\
\ds\ \ = \!\! \sum_{k=0}^{\min\{a,b\}}\!\!\!\!
\frac{1}{k!(a-k)!(b-k)!}\
\overline{{\rm Sym}}^{a}_{t}\overline{{\rm Sym}}^{b}_{s}
\left(
P\sk{f_\a(t_{1})\cdots  f_\a(t_{a-k})
f_{\a+\b}(t_{a-k+1})\cdots f_{\a+\b}(t_{a})}\right.
\\
\ds\qquad \qquad \times\ P\sk{f_\b(s_{k+1})\cdots
f_\b(s_{b})}
\left.Y(q^{-1}t_{a-k+1},\ldots,q^{-1}t_a;s_1,\ldots,s_k)
\right).
\end{array}
\ee
\end{theor}
This theorem reduces calculation the weight function to
calculating the projections of strings.

\subsection{Projections of strings}
We first describe projections of single currents.
 For any current $a(z)=\sum_{n\in\ZZ} a_nz^{-n}$, let
$a^\pm(z)$ denote the currents ($a(z)=a^+(z)-a^-(z)$)
\bn\label{pro-int}
\begin{array}{l}
a^+(z)=\oint \frac{a(w)}{1-\frac{w}{z}}\frac{dw}{z}\ = \
\sum_{n>0}a_nz^{-n} ,\\
a^-(z)=-\oint \frac{a(w)}{1-\frac{z}{w}}\frac{dw}{w}\ = \
-\sum_{n\leq 0}a_nz^{-n} .
\end{array}
\ee
\begin{prop}\label{single}
Projections of the currents $f_\a(z)$, $f_\b(z)$, and $f_{\a+\b}(z)$
can be written as
\be
\begin{array}{l}
P\sk{f_{\a}(t)}=f^+_\a(t)\, ,\qquad
P\sk{f_{\b}(t)}=f^+_\b(t)\,,\\
\ds P\sk{f_{\a+\b}(t)}=S_{\b}(f^+_\a(t))\ =\
f^+_\a(t)f_\b[0]-qf_\b[0]f^+_\a(t)\,.
\end{array}
\ee
\end{prop}
There are also analogous formulas for the opposite projection:
\be\label{3.14}
\begin{array}{l}
P^-\sk{f_{\a}(t)}=-f^-_\a(t)\,,\qquad P^-\sk{f_{\b}(t)}=-f^-_\b(t)\,, \\
\ds P^-\sk{f_{\a+\b}(t)}=q^{-1}S_\a(f_\b^-(q^{-1}t))= q^{-1}f^-_\b(q^{-1}t)
f_\a[0]-f_\a[0]f^-_\b(q^{-1}t)\,.
\end{array}
\ed
We define a set   of  rational functions
of the variable $s$ depending on parameters  $s_1,\ldots,s_b$:
\be\label{RF15}
\varphi_{s_j}(s;s_1,\ldots,s_{b})=\prod_{i=1,\ i\neq j}^{b}
\frac{s-s_i}{s_j-s_i}\prod_{i=1}^{b}\frac{q^{-1}s_j-qs_i}{
q^{-1}s-qs_i}\ .
\ee
As functions of $s$, they have simple poles at the points $s=q^2s_i$,
$i=1,\ldots,b$, tend to zero as  $s\to\infty$,
and have properties $\varphi_{s_j}(s_i;s_1,\ldots,s_b)=\delta_{ij}$.
Set  \r{RF15} is uniquely defined by these properties.

We define the combination of currents
\be\label{long-cur}
f_\gamma(t;t_{1},\ldots,t_b)=f_\gamma(t)-
\sum_{m=1}^b \varphi_{t_m}(t;t_{1},\ldots,t_b)f_\gamma(t_m) \,,
\ee
where $\gamma$ coincides either with the simple root $\alpha$, the simple root
$\beta$, or the composite root $\alpha+\beta$.

\begin{theor}\label{prstring}
The projection of the string \rf{string} has the factored form
\be\label{po33}
\begin{array}{l}
\ds P\sk{f_\a(t_1)\cdots f_\a(t_{a-k})f_{\a+\b}(t_{a-k+1})\cdots
f_{\a+\b}(t_{a})}=
\\
\ds\quad = \prod_{1\leq i\leq a-k <j\leq a}
\frac{qt_i-q^{-1}t_j}{t_i-t_j}
\prod_{1\leq i<j\leq a}
\frac{q^{-1}t_i-qt_j}{qt_i-q^{-1}t_j}\\
\ds\qquad\times\
P\sk{f_{\a+\b}(t_a)}P\sk{f_{\a+\b}(t_{a-1};t_a)}
\cdots P\sk{f_{\a+\b}(t_{a-k+1};t_{a-k+2},\ldots,t_a)}\\
\ds\qquad\quad\times\
P\sk{f_{\a}(t_{a-k};t_{a-k+1},\ldots,t_a)}\cdots P\sk{f_{\a}(t_{1};t_{2},\ldots,t_a)}.
\end{array}
\ee
\end{theor}

\subsection{Examples}\label{sect-ex}
We give several explicit examples illustrating Theorems~\ref{act-cal}
 and \ref{prstring}.
The second, third and forth examples are given with the
corollary to Theorem~\ref{sl2prop} taken into account:
$$
\begin{array}{rcl}
\ds{P}\sk{f_\a(t_1)f_{\a+\b}(t_2)}&=&\ds\frac{q^{-1}t_1-qt_2}{t_1-t_2}\
P\sk{f_{\a+\b}(t_2)}\sk{f_\a^+(t_1)-\frac{(q-q^{-1})t_2}{qt_2-q^{-1}t_1}
f^+_\a(t_2)}\, ,\\[5mm]
\ds{P}\sk{f_\a(t_1)f_{\a}(t_2)}&=&\ds
f_{\a}^+(t_1)\sk{f^+_\a(t_2)-\frac{(q-q^{-1})t_1}{qt_1-q^{-1}t_2}
f^+_\a(t_1)}\, ,\\[5mm]
\ds{P}\sk{f_{\a+\b}(t_1)f_{\a+\b}(t_2)}&=&
\ds {P}\sk{f_{\a+\b}(t_1)}\sk{{P}\sk{f_{\a+\b}(t_2)}-\frac{(q-q^{-1})t_1}{qt_1-q^{-1}t_2}
{P}\sk{f_{\a+\b}(t_1)}}\, ,
\end{array}
$$
$$
\begin{array}{rcl}
\ds{P}\sk{f_\a(t_1)f_{\a}(t_2)f_{\a}(t_3)}&=&\ds
f_{\a}^+(t_1)\sk{f^+_\a(t_2)-\frac{(q-q^{-1})t_1}{qt_1-q^{-1}t_2}
f^+_\a(t_1)}\times\\[5mm]
&\times&\ds \left(f_\a^+(t_3)
-\frac{t_1-t_3}{t_2-t_3}\frac{(qt_1-q^{-1}t_3)(q-q^{-1})t_3}
{(qt_1-q^{-1}t_3)(qt_2-q^{-1}t_3)}f_\a^+(t_2)-\right.\\[5mm]
&-&\left.\ds\frac{t_2-t_3}{t_2-t_1}\frac{(qt_3-q^{-1}t_1)(q-q^{-1})t_1}
{(qt_1-q^{-1}t_3)(qt_2-q^{-1}t_3)}f_\a^+(t_1)\right)
\end{array}
$$
and
\be\label{exam2}
\begin{array}{l}
\ds P\sk{f_\a(t_1)f_\a(t_2)f_\b(s_1)f_\b(s_2)}=
P\sk{f_\a(t_1)f_\a(t_2)} P\sk{f_\b(s_1)f_\b(s_2)}+\\
\ds\quad +\ \overline{\rm Sym}_{s_1,s_2}\sk{
\overline{\rm Sym}_{t_1,t_2}
\sk{{P}\sk{f_\a(t_1)f_{\a+\b}(t_2)}\frac{t_2}{t_2-qs_1}}f^+_\b(s_2)}+\\
\ds\quad +\ \frac{1}{2}\
\overline{\rm Sym}_{t_1,t_2}\sk{
{P}\sk{f_{\a+\b}(t_1)f_{\a+\b}(t_2)}
\overline{\rm Sym}_{s_1,s_2}\sk{\frac{t_1}{t_1-qs_1}
\frac{t_2}{t_2-qs_2}\frac{qt_2-s_1}{t_2-qs_1}}}.
\end{array}
\ee
We note that in Theorem~\ref{prstring} and in the examples considered above,
the normal ordering of the
roots is changed from $\a,\a+\b,\b$ to $\a+\b,\a,\b$. The correct
 normal ordering $\a,\a+\b,\b$ can be restored by using the commutation
 relations given in Proposition~\ref{restore-ord}.
The result of this calculations is that the second line in the formula
 \r{exam2} can be  replaced with the expression
$$ \ds\quad \ \overline{\rm Sym}_{s_1,s_2}\sk{
\overline{\rm Sym}_{t_1,t_2}
\sk{f^+_\a(t_1)P({f}_{\a+\b}(t_2;t_1))\ \frac{q^{-1}t_1-qt_2}{t_1-t_2}\
\frac{qt_1-s_1}{t_1-qs_1}\
\frac{t_2}{t_2-qs_1}
}f^+_\b(s_2)}. $$

\subsection{Universal weight function for $\Uqdva$}

The currents $e_\a(z)$, $f_\a(z)$, and $\psi^\pm_\a(z)$, as well as Chevalley
generators $e_{\pm \a_i}$ and $k_{\a_i}^\pm$, $i=0,1$,  generate a Hopf subalgebra
$\Uqdva$ in $\Uqtri$. For this algebra, the weight function and projection operators
can be defined independently. We observe that the corresponding
projection of the product $f_\a(t_1)\cdots f_\a(t_n)$ coincides with its
projection inside the algebra $\Uqtri$. As a corollary of Theorem~\ref{prstring}, we
obtain the description of the weight function for $\Uqdva$. It also admits a
simple integral presentation.
\begin{theor}\label{sl2prop}${}$
\begin{itemize}
\item[(i)]
Universal weight function \rf{Wts} can be written as
\be\label{po333}
\begin{array}{rcl}
\ds W(t_1,\ldots, t_a)&=&
\ds P\sk{f_\a(t_1)\cdots f_\a(t_{a})}\ =\\
&=&
\prod_{1\leq i<j\leq a}
\frac{q^{-1}t_i-qt_j}{qt_i-q^{-1}t_j}\
%\prod^{\longleftarrow}_{a\geq \ell \geq 1}
f^+_\a(t_a)f^+_\a(t_{a-1};t_a)\cdots
f^+_\a(t_1;t_{2},\ldots,t_a)\,.
\end{array}
\ee
\item[(ii)]
Weight function \rf{po333} admits the integral representation
\be\label{intt1}
\!P\sk{f_\a(t_1)\cdots f_\a(t_a)}=
\prod_{i<j}\frac{t_i-t_j}{qt_i-q^{-1}t_j}
\mathop{\oint\!\cdots\!\oint}%\limits_{ w_a\gg \cdots \gg w_1}
Z(\overline{t};\,\overline{w})\,f_\a(w_a)\frac{dw_a}{w_a}\cdots f_\a(w_1)
\frac{dw_1}{w_1}\,.
\ed
\end{itemize}
\end{theor}
The currents $f^+_\a(t_\ell;t_{\ell+1},\ldots,t_a)$
 are defined by formula  \r{long-cur}, and
the   kernel
$Z(\overline{t};\,\overline{w})$ of the integral transformation is defined in
\r{rat-Y}.

{\it Proof.}\ Statement (i) is a particular case of
Theorem~\ref{prstring}. Statement (ii) is obtained from (i) by substituting
expressions \rf{pro-int} and the elementary identity
$$
\frac{1}{t-w}-\sum_{m=1}^b\varphi_{t_m}(t;t_1,\ldots,t_b)\ \frac{1}{t_m-w}
=\frac{1}{t-w}\prod_{i=1}^b\frac{t-t_i}{w-t_i}\ \frac{q^{-1}w-qt_i}{q^{-1}t-qt_i}
\ .
$$

We note that because the factor before the integral in  \rf{intt1} has the same
$q$-symmetry properties as the product of currents
$f_\a(t_1)f_\a(t_2)\cdots f_\a(t_a)$,
the integral is itself symmetric under permutations
of the parameters $t_1,\ldots,t_a$. This means that we can use
the kernel $$Z(t_{\sigma(1)},\ldots,t_{\sigma(a)};w_1,\ldots,w_a)$$ in \rf{intt1}
instead of the kernel $Z(t_1,\ldots,t_a;w_1,\ldots,w_a)$
for any permutation $\sigma\in S_a$.

\begin{cor}\label{inverse}
The projection of the product of currents can be written
in the ``direct'' order
\be\label{intt3}
\begin{array}{rcl}
\!\!\!P\sk{f_\a(t_1)\cdots f_\a(t_a)}\!&=&
\!\!\prod_{1\leq i<j\leq a}
\frac{t_i-t_j}{qt_i-q^{-1}t_j}
\mathop{\oint\!\cdots\!\oint}
Z(^\op\overline{t};\,^\op\overline{w})\,f_\a(w_1)\frac{dw_1}{w_1}\cdots f_\a(w_a)
\frac{dw_a}{w_a}\\
&=&f^+_\a(t_1)f^+_\a(t_2;t_1)\cdots
f^+_\a(t_{a-1};t_1,\ldots,t_{a-2})f^+_\a(t_a;t_1,\ldots,t_{a-1}).
\end{array}
\ee
\end{cor}
To prove this corollary, it is sufficient to rename the parameters in integral
 \r{intt1}
 $t_i\to t_{a+1-i}$, $i=1,\ldots,a$ and calculate the integral, or to rename
 the variables in \rf{po333} and analytically continue the result
 to the original domain.

\subsection{Combinatorial  identity for the kernels $Y(\overline{t};\overline{s})$ and
 $Z(\overline{t};\overline{s})$}

The opposite current Borel subalgebra $U_E$ (see Sec.~\ref{bor-sub})
also admits a decomposition into a product of its intersections with
Borel subalgebras,
$$%%\be\label{F-dec}
U_e^+=U_E\cap U_q(\mathfrak{b}_+)=U_E\cap U_q(\mathfrak{n}_+)\,,\qquad
U_E^-=U_E\cap U_q(\mathfrak{b}_-)\ ,
$$%%\ee
such that the relations
\bn \label{Pdef2}
 \tilde{P}^+(e_1e_2)=e_1\varepsilon(e_2)\,,
\qquad \tilde{P}(e_1e_2)= \
\tilde{P}^-(e_1e_2)=\varepsilon(e_1)e_2\, , \ed where  $e_1\in
U_e^+$ and $e_2\in U_E^-$, define projection operators
$\tilde{P}^\pm$ that  map $U_E$ to their subalgebras $U_E^-$ and
$U_e^+$ and have properties analogous to the properties of the
projection operator $P^\pm$.

In particular, the projection $\tilde{P}\sk{e_\a(s_1)\cdots e_\a(s_b)}$
 admits an integral representation with the factored kernel
 $Z(\overline{w};\,\overline{t})$ (see \r{rat-Y}):
\be\label{intt4}
\!\tilde{P}\sk{e_\a(s_1)\cdots e_\a(s_a)}=
\prod_{i<j}\frac{s_i-s_j}{q^{-1}s_i-q^{}s_j}
\mathop{\oint\!\cdots\!\oint}%\limits_{ w_a\gg \cdots \gg w_1}
{Z}(\overline{w};\overline{s})e_\a(w_1)\frac{dw_1}{w_1}\cdots e_\a(w_a)
\frac{dw_a}{w_a}\, .
\ed

We let the symbol  $\widetilde{\rm Sym}_{s}^{n}g(s_1,\ldots,s_n)$
denote $q^{-1}$-symmetrization of a function $g(s_1,\ldots,s_n)$:
\be
\begin{array}{rcl}\label{qs3r}
\widetilde{\rm Sym}_{s}^{n}\ g(s_1,\ldots,s_n) &=&
\sum_{\nu\in S_{n}}
\prod_{\ell<\ell'\atop \nu^{-1}(\ell)>\nu^{-1}(\ell')}
\frac{q^{}s_\ell-q^{-1}s_{\ell'}}
{q^{-1}s_\ell-q^{}s_{\ell'}}\
g(s_{\nu(1)},\ldots, s_{\nu(n)})=\\&=&
\prod_{i<j}\frac{qs_i-q^{-1}s_j}{q^{-1}s_i-q^{}s_j}\
\overline{\rm Sym}_{s}^{n}\ g(s_n,\ldots,s_1)\,.
\end{array}
\ee

The current Borel subalgebras $U_F$ and $U_E$ are Hopf dual with respect to the
Hopf structure $\Delta^{(D)}$. Let $\<\ ,\,\>:U_E\otimes U_F\to\CC$
be the corresponding Hopf pairing. It has the properties
 $\<ab,x\>=\<b\otimes a,\Delta^{(D)}(x)\>$ and
$\<a,xy\>=\<\Delta^{(D)}(a),x\otimes y\>$; for the algebra  $\Uqdva$, it is given by
\bn
%\begin{array}{rcl}
\<e_\a(s_1)\cdots e_\a(s_n),f_\a(t_1)\cdots f_\a(t_n)\>=%&=&
(q^{-1}-q)^{-n}\  \overline{\rm Sym}^n_t \sk{\prod_{k=1}^n
\delta\left(\frac{s_k}{t_k}\right)}. %=\\
%&=&
%(q^{-1}-q)^{-n}\  \widetilde{\rm Sym}^n_s \sk{\prod_{k=1}^n
%\delta\left(\frac{t_k}{s_k}\right)}.
%\end{array}
\label{pr}
\ed
We note that for obvious reasons, the right-hand side of equality
 \r{pr} can be rewritten as  a $q^{-1}$-symmetri\-za\-tion in the variables
 $s_1,\ldots,s_n$ in the domain
$|s_1|\ll|s_2|\ll\ldots\ll|s_n|$.
\begin{prop} \label{adjoint}
The operators $P^\pm$ and $\tilde{P}^\mp$ are adjoint with respect to the Hopf
pairing $\<\ ,\ \>$: for any $f\in U_F$ and $e\in U_E$ we have
$$
\<e,P^+(f)\>=\<\tilde{P}^-(e),f\>\,, \qquad
\<e,P^-(f)\>=\<\tilde{P}^+(e),f\>\, .
$$
\end{prop}

{\it Proof}.\
We let
$\tilde\mathcal{R}\in U_E\otimes U_F$ denote the tensor of Hopf pairing \r{pr}.
In the notation in \cite{DKhP},  $\tilde\mathcal{R}$ coincides with
 $\sk{\mathcal{R}^{21}}^{-1}$.
It was established in Sec. 4.2 \
 in that paper that the two pairs
$(U_f^-,U_F^+)$ and $(U_e^+,U_E^-)$ of subalgebras of current Borel
algebras form a biorthogonal decomposition of a quantum affine
algebra (see Sec. 4.1 in \cite{DKhP},  for the definition). This
implies that the tensor $\tilde{\cal R}$ $\in$ $U_E\ot U_F$ of the
Hopf pairing admits a decomposition
 $\tilde{\cal R}$ $=$ ${\cal R}_1 {\cal R}_2$, where
$${\cal R}_1= (1\ot P^-)\tilde{\cal R}=(\tilde{P}^+\ot 1)\tilde{\cal R}\,,\quad
{\cal R}_2= (1\ot P^+)\tilde{\cal R}=(\tilde{P}^-\ot 1)\tilde{\cal R}\, ,$$
such that $e\in U_E$, $f\in U_F$ and the equalities
\bn\label{Pad}
\<e,P^+(f)\>=\<\tilde{P}^-(e),f\>=\<{\cal R}_2, f\ot e\>\,, \quad
\<e,P^-(f)\>=\<\tilde{P}^+(e),f\>=\<{\cal R}_1, f\ot e\>\,
\ed
hold.
\hfill{$\Box$}

\begin{prop}${}$\label{PPprop}

\noindent  (i) For any sets of variables $\overline{t}=\{t_1,...,t_n\}$
 and $\overline{s}=\{s_1,...,s_n\}$, we have the equalities
%Спаривание $ \<P(f_\a(t_1)\cdots f_\a(t_n)),
%e_\a(s_1)\cdots e_\a(s_n)\>$ равно:
\begin{equation}\label{is-pair1}
%\begin{array}{c}
 \<e_\a(s_1)\cdots e_\a(s_n),P(f_\a(t_1)\cdots f_\a(t_n))\>=
%\gamma
(q^{-1}-q)^{-n}\prod_{i<j}\frac{t_i-t_j}{qt_i-q^{-1}t_j}
\widetilde{\rm Sym}_{s}^{n}Z(\overline{t};\,^\op\overline{s})\, ,
\end{equation}
\begin{equation}\label{is-pair2}
 \<\tilde{P}(e_\a(s_1)\cdots e_\a(s_n)),f_\a(t_1)\cdots f_\a(t_n)\>=
%\gamma
(q^{-1}-q)^{-n}\prod_{i<j}\frac{s_i-s_j}{q^{-1}s_i-q^{}s_j}
\overline{\rm Sym}_{t}^{n}Z(\overline{t};\overline{s})\,.
%\end{array}
\end{equation}
\noindent  (ii)
For any two permutations $\sigma,\sigma'\in S_n$,
we have the following identity in the ring of  functions, symmetric with respect to both sets of variables
$\overline{t}$ and $\overline{s}$:
\begin{equation}
%\begin{eqnarray}\label{idZZ1}%{l}
%\ds
\prod_{i<j}\frac{qt_i-q^{-1}t_j}{t_i-t_j}\
\overline{\rm Sym}_{t}^{n}\ Z(\overline{t};\,^\sigma\overline{s})%\\
%\qquad
=%\ds
\prod_{i<j}\frac{qs_i-q^{-1}s_j}{s_i-s_j}\
\overline{\rm Sym}_{s}^{n}\ Z(^{\sigma'}\overline{t};\overline{s})\, ,%%\\
\end{equation}
\begin{equation}\label{idZZ10}
%\ds
\prod_{i<j}\frac{qt_i-q^{-1}t_j}{t_i-t_j}\
\overline{\rm Sym}_{t}^{n}\ Y(\overline{t};\,^\sigma\overline{s})%\\
%\qquad
=%\ds
\prod_{i<j}\frac{qs_i-q^{-1}s_j}{s_i-s_j}\
\overline{\rm Sym}_{s}^{n}\ Y(^{\sigma'}\overline{t};\overline{s})\, .
%\end{eqnarray}
\end{equation}
\end{prop}

{\it Proof.}\
 Statement (i) can be obtained by substituting integral
presentations \rf{intt1} and  \rf{intt4} in the corresponding Hopf pairings.
From (i), after the $q^{-1}$- symmetrization is replaced with
$q$-symmetrization according to \rf{qs3r},  we obtain (ii) for
 $\sigma = \sigma' =1$
(we recall that the functions
$Z(\overline{t};\overline{s})$ and $Y(\overline{t};\overline{s})$
differ by a simple factor, symmetric with respect to both sete of variables; see
 \r{rat-Y}). Next,
both sides of equality \rf{is-pair1} are $q$-symmetric with respect
to the variables $t$. Because the product
$\prod_{i<j}\frac{t_i-t_j}{qt_i-q^{-1}t_j}$ is also $q$-symmetric,
the remaining factor $\widetilde{\rm
Sym}_{s}^{n}Z(\overline{t},^\op\overline{s})$ is symmetric with
respect to the variables $t$, and for any$\sigma\in S_n$, we hence
have
\begin{equation}
\overline{\rm Sym}_{s}^{n}Z(\overline{t};\overline{s})=
\overline{\rm Sym}_{s}^{n}Z(^\sigma\overline{t};\overline{s})\quad
\quad
\overline{\rm Sym}_{t}^{n}Z(\overline{t};\overline{s})=
\overline{\rm Sym}_{t}^{n}Z(\overline{t};\,^\sigma\overline{s})\ ,
\end{equation}
which implies the statement (ii).
\hfill{$\Box$}

Identity \r{idZZ1}
has several  proofs via direct calculations. Our proof is based on interpreting
both sides of this identity as specific matrix elements of a $\Uqdva$
weight function. In what follows, we use identity \rf{idZZ10} in the form
\begin{equation}\label{idZZ1}
%\ds
\prod_{i<j}\frac{qt_i-q^{-1}t_j}{t_i-t_j}\
\overline{\rm Sym}_{t}^{n}\ Y(t_1,...,t_n;s_n,...,s_1)%\\
%\qquad
=%\ds
\prod_{i<j}\frac{qs_i-q^{-1}s_j}{s_i-s_j}\
\overline{\rm Sym}_{s}^{n}\ Y(t_1,...,t_n;s_1,...,s_n)\,.
%\end{eqnarray}
\end{equation}

We let $\mathcal{Z}(\overline{t},\overline{s})$ denote the left- or the right-hand
side of identity  \r{idZZ1} divided by the product
$\prod_{i=1}^n t_i$. As follows from the above considerations,
this function $\mathcal{Z}(\overline{t},\overline{s})$ is symmetric with respect to both sets
of variables $\overline{t}$ and $\overline{s}$ and has a ``physical'' meaning.
It coincides with the partition function of the complete inhomogeneous six-vertex model
on a square lattice with domain-wall boundary conditions. As pointed out to us by N.~Slavnov,
this function has the determinant representation
\be\label{stat-s}
\mathcal{Z}(\overline{t},\overline{s})=\frac{\ds \prod_{i,j=1}^n(qt_i-q^{-1}s_j)}
{\ds \prod_{i<j}(t_i-t_j)(s_j-s_i)}
\det\left|\frac{1}{(t_i-s_j)(qt_i-q^{-1}s_j)}\right|_{i,j=1,\ldots,n}\ .
\ee

\setcounter{equation}{0}

\section{Analytic properties of strings}\label{anal-st}

\subsection{Properties of the current $f_{\a+\b}(z)$}
\label{compcurrent}

The current $f_{\a+\b}(z)$ is defined in Sec.~\ref{com-cs} by
relation
 \rf{rec+}. We first note that because of relations \rf{1}, it admits the analytic representation
\be\label{c-c-f}
f_{\a+\b}(z)\ = \ \res{w=z}f_{\a}(w)f_{\b}(q^{-1}z)\,\frac{dw}{w}\
=\ (q-q^{-1})
f_\b(q^{-1}z)f_\a(z)\,
\ee
in addition to \rf{rec+}. We can also obtain the current $f_{\a+\b}(z)$ as a result of
the adjoint action related to the comultiplication $\D^{(D)}$.
We define the left and right adjoint actions with respect to the coalgebra structure
$\Delta^{(D)}$:
\bn\label{adcurrent}
{\rm ad}_x^{(D)}(y)=\sum_ja^{(D)}(x'_j)\cdot y\cdot  x''_j\,,
\qquad
\tilde{\rm ad}_x^{(D)}(y)=\sum_j  x''_j\cdot y\cdot(a^{(D)})^{-1}(x'_j)
\ed
if
$\Delta^{(D)}(x)=\sum_jx'_j\otimes x''_j$.
We call them the {\it current adjoint actions}. We have
\bn
\begin{array}{l}
{\rm ad}^{(D)}_{f_{i}(z)}(y)=yf_{i}(z)-f_{i}(z)\sk{\psi^+_i(z)}^{-1}y\,
\psi^+_i(z)\, ,\\
\tilde{\rm ad}^{(D)}_{f_{i}(z)}(y)=f_{i}(z)y-\psi^+_i(z)\,y\,
\psi^+_i(z)^{-1}f_{i}(z)\,.
\end{array}
\ed
\begin{prop}${}$
 We have the equalities
\bn\label{adac4}
{\rm ad}^{(D)}_{f_{\b}(w)}\sk{f_{\a}(z)}=\d(zq^{-1}/w)f_{\a+\b}(z)\,,\qquad
\tilde{\rm ad}^{(D)}_{f_{\a}(w)}\sk{f_{\b}(z)}=\d(qz/w)f_{\a+\b}(qz)\,,
\ed
and hence
\be\label{adac1}
f_{\a+\b}(z)=\oint \frac{dw}{w}\ {\rm ad}^{(D)}_{f_{\b}(w)}\sk{f_{\a}(z)}\ = \
\oint \frac{dw}{w}\ \tilde{\rm ad}^{(D)}_{f_{\a}(w)}\sk{f_{\b}(q^{-1}z)}\,.
\ee
\end{prop}

\begin{prop} \label{prop4.2}
\cite{DKh} The following relations hold in $\UqtriD:$
\be
\begin{array}{c}
\label{com-aa}
f_{\a}(z)f_{\a+\b}(w)\ =\  \frac{q^{-1}z-qw}{z-w}f_{\a+\b}(w)f_{\a}(z)\, ,\\
 f_{\a+\b}(qz)f_{\b}(w)\ =\
\frac{q^{-1}z-qw}{z-w}f_{\b}(w)f_{\a+\b}(qz)\, ,\\
 \frac{qz-q^{-1}w}{z-w}f_{\a+\b}(z)f_{\a+\b}(w)\ =\
 \frac{zq^{-1}- qw}{z-w}f_{\a+\b}(w)f_{\a+\b}(z)\,.
\end{array}
\ed
\end{prop}
We note that in the analytic language, both sides of all relations
\rf{com-aa} are analytic functions in $\sk{\CC^*}^2$. This means,
for instance, that the product $f_{\a}(z)f_{\a+\b}(w)$ has no zeroes
and poles, while the product $f_{\a+\b}(w)f_{\a}(z)$ has a simple
zero at $z=w$ and a simple pole at $z=q^2w$.

\noindent {\it Proof.} The proof combines relations \rf{1} and
the Serre relations in the analytic form \cite{E}.
Namely, in the algebra $\UqtriD$,
\begin{itemize}
\item[(i)] the products
$f_\a(z)f_\a(w)$ and $f_\b(z)f_\b(w)$ have a simple zero at  $z=w$, and
\item[(ii)]  the products
$(z_1-qz_2)(z_2-qz_3)(z_1-q^{-2}z_3)f_\a(z_1)f_\b(z_2)f_\a(z_3)$
and
$(z_1-qz_2)(z_2-qz_3)(z_1-q^{-2}z_3)f_\b(z_1)f_\a(z_2)f_\b(z_3)$
vanish on the lines $z_2=qz_1=q^{-1}z_3$ and $z_2=q^{-1}z_1=q^{}z_3$.
\end{itemize}

The properties of products of currents given in
Proposition \ref{prop4.2} admit a straightforward generalization
 to strings.
\begin{prop}\label{hyper}${}$

\noindent (i)
Strings \rf{string}
% $f_{\a+\b}(u_{a+b})\cdots \!f_{\a+1}(u_{a+b})f_\a(u_a)$
%and $f_\b(v_{a+b})\cdots \!f_\b(v_{a+1})
%f_{\a+\b}(v_{a+1})\cdots \!f_{\a+\b}(v_{1})$
only have simple poles on the hyperplanes
 $u_i=q^{-2}u_j$ and simple zeros on the hyperplanes $u_i=u_j$,
 where either  $1\leq i<j\leq a$ or $a+1\leq i<j\leq a+b$.

\noindent(ii) Opposite strings \rf{stringop}
only have simple poles on the hyperplanes
 $u_i=q^{2}u_j$ and simple zeros on the hyperplanes $u_i=u_j$ for all pairs
 $i<j$.

\noindent (iii) The strings and opposite strings are related.
In particular,
\be
\begin{array}{c}\label{oop}
f_{\a}(t_{k+1})\cdots f_{\a}(t_a)f_{\a+\b}(t_1)\cdots\! f_{\a+\b}(t_k)=\\
=\!\prod_{i,j:\ 1\leq i\leq k<j\leq n}\frac{q^{-1}t_i-qt_j}{t_i-t_j}
f_{\a+\b}(t_1)\cdots\! f_{\a+\b}(t_k)
f_{\a}(t_{k+1})\cdots\! f_{\a}(t_n)\,.
\end{array}
\ee
\end{prop}

\subsection{Screening operators and projections of
$f_{\a+\b}(z)$}
Let $\tilde{S}_i$ denote the screening operators
$\tilde{S}_i=\tilde{\rm ad}_{\fin{i}{0}}$ with respect to the adjoint action $\tilde{\rm ad}$
 in $\UqtriD$:
 $\tilde{\rm ad}_x(y)=\sum_ix''_i\cdot y\cdot a(x'_i)$,
where
$\Delta(x)=\sum_ix'_i\otimes x''_i$, such that
\bn\label{coad1}
\ds \tilde S_{i}\sk{y}= \ds
\fin{i}{0}\,y- k^{-1}_{i}\,y\,k_{i}\,\fin{i}{0}\,.
\ed
 We have the following relations between screening and projection operators.
\begin{prop}${}$\label{PS}
\begin{itemize}\label{zero-mode-com}
\item[(i)] For any $y\in U_F$ and $i=\a,\b$, we have equalities
\be P^+\sk{\oint\frac{dw}{w}\ {\rm ad}^{(D)}_{f_{i}(w)}(y)}=S_i\sk{P^+(y)}\,,\quad
P^-\sk{\oint\frac{dw}{w}\ \tilde{\rm ad}^{(D)}_{f_{i}(w)}(y)}=\tilde{S}_i\sk{P^-(y)}\,.
\ed
\item[(ii)]
\ The screenings operators $ S_{i}$ and $\tilde S_{i}$ are related as
\be\label{ad-rel}
\tilde S_{i}\sk{y}= -q^{-2}k^{-1}_{i}\,S_{i}
\sk{y}\,k_{i}\,.
\ee
\item[(iii)]\ The screening operators $S_{i}$ and
 $\tilde S_{i}$ commute with  projections $\Pfpm$: for any $y\in
U_F$,
\be\label{ad-com}
\Pfpm S_{i}\sk{y}=S_{i}\Pfpm\sk{y}\,,\quad
\Pfpm \tilde S_{i}\sk{y}=\tilde S_{i}\Pfpm\sk{y}\,.
\ee
\end{itemize}
\end{prop}
\noindent
{\it Proof.}\ \ Statements (i) and (ii) are obvious.
 We prove the equality $P S_{f_i[0]}\sk{y}=S_{f_i[0]}P\sk{y}$, where
 $y=y_1y_2$ and
$y_1\in U_f^-$ and $y_2\in U_F^+$. The adjoint action has the property
${\rm ad}_x(y_1\cdot y_2)=\sum_i {\rm ad}_{x'_i}(y_1)\cdot {\rm ad}_{x''_i}(y_2)$, if
$\Delta(x)=\sum_ix'_i\otimes x''_i$. This
implies the relation
\bn
 S_{i}\sk{y}= S_{i}\sk{y_1}k_iy_2k_i^{-1}+ y_1  S_{i}\sk{y_2},
\ed
$$PS_{i}\sk{y}\ = \
\ve\sk{S_{i}\sk{y_1}}k_iy_2k_i^{-1}+\ve(y_1)S_{i}\sk{y_2}=S_{i}P\sk{y}\, $$
because the screening operator $S_{i}$ preserves the subalgebras
$U_F^+$ and $U_f^-$ and $\ve( S_{i}(y))=0$ for any $y\in\U_f$ except $y=1$.
\hfill$\Box$
\medskip

The properties of the screenings operators and of the adjoint
actions allow calculating the projections of the current
$f_{\a+\b}(z)$ and establishing the corresponding normally ordered
decomposition \rf{addF} for it.
\begin{prop}${}$
\begin{itemize}
\item[(i)]
The projections of the current $f_{\a+\b}(z)$ are
\be\label{iter2}
\Pfp\sk{f_{\a+\b}(z)}=
S_{\b}\sk{f^+_\a(z)}\,,\quad
\Pfm\sk{f_{\a+\b}(z)}=-
\tilde S_{\a}\sk{f^-_{\b}(q^{-1}z)}\,.
\ee
\item[(ii)] We have the normally ordered expansion
\be\label{pl-r2}
 f_{\a+\b}(z) -P\sk{f_{\a+\b}(z})=
(q^{-1}-q)\sk{f^-_{\b}(q^{-1}z)f_{\a}(z)}^+ - f^-_{\a+\b}(z)\,.
\ee
\end{itemize}
\end{prop}
{\it Proof.}\ \ Statement (i) follows from Proposition \ref{PS} and relation
\rf{adac1}.
Next, formal integral
 \r{rec-f} can be written as %\rf{comf-m},
\bn\begin{array}{rcl}\label{comf-m}
\ds f_{\a+\b}(w)&=&
f_{\a}(w)\fin{\b}{0}-q^{-1}\fin{\b}{0}f_{\a}(w)
-(q^{-1}-q)\sum_{k>0}\fin{\b}{-k}f_{\a}(w)\sk{q^{-1}w}^k=\\
&=&f_{\a}(w)\fin{\b}{0}-q^{}\fin{\b}{0}f_{\a}(w)
-(q^{-1}-q)\sum_{k\geq 0}\fin{\b}{-k}f_{\a}(w)\sk{q^{-1}w}^k= \\
&=&S_{\b}\sk{f_{\a}(w)}
+(q^{-1}-q)f^-_{\b}(q^{-1}w)f_{\a}(w)\,.
\end{array}
\ed
Applying the operation $\oint(\cdot)\frac{dw}{z(1-w/z)}$ see
\rf{pro-int} to both sides of \rf{comf-m} proves Statement (ii).
\hfill$\Box$

\subsection{Proof of Theorem~\ref{prstring}}

We first calculate the projection of the opposite string $P(f_{\a+\b}(t_1)\cdots f_{\a+\b}(t_k)
f_{\a}(t_{k+1})\cdots f_{\a}(t_a))$.

\begin{prop}${}$\label{prop4.6}

\noindent (i)
The projection of the string $f_{\a+\b}(t_1)\cdots \!f_{\a+\b}(t_k)
f_{\a}(t_{k+1})\cdots \!f_{\a}(t_a)$ with $a>k$
 can be written as
\be\label{RF11}
\begin{array}{l}
\ds {P}\sk{f_{\a+\b}(t_1)\cdots f_{\a+\b}(t_k)
f_{\a}(t_{k+1})\cdots f_{\a}(t_a)}=\\
\ds\quad  = {P}\sk{f_{\a+\b}(t_1)\cdots f_{\a+\b}(t_k)
f_{\a}(t_{k+1})\cdots f_{\a}(t_{a-1})}
f^+_\a(t_a)+
\sum_{j=1}^{a-1}\frac{X_j(t_1,\ldots,t_{a-1})} {t_a-q^2t_j}\,.
\end{array}
\ee
%where operators $X_j$ do not depend on the variable $t_a$;
\noindent (ii)
The projection of the string $f_{\a+\b}(t_1)\cdots f_{\a+\b}(t_k)$
can be represented as
\be\label{RF12}
\begin{array}{l}
\ds {P}\sk{f_{\a+\b}(t_1)\cdots \! f_{\a+\b}(t_k)}
\ds  = {P}\sk{f_{\a+\b}(t_1)\cdots \!f_{\a+\b}(t_{k-1})}P(f_{\a+\b}(t_{k}))+
\sum_{j=1}^{k-1}\frac{X'_j(t_1,...,t_{k-1})} {t_k-q^2t_j}\,.
\end{array}
\ee
%where operators ${X'}_j$ do not depend on the variable $t_k$;
%\end{itemize}
\end{prop}
{\it Proof.}\ Relation \rf{RF11} is proved based on inductively
using the following lemma, which shows that as the current
$f_\a^-(t_a)$ is moved to the left of the string, only simple poles
at the points $t_a=q^2t_j$ appear and the corresponding
operator-valued coefficient $X_j$ at $(t_a-q^2t_j)^{-1}$  is
therefore independ of $t_a$.
\begin{lem}\label{co-min} The following relations hold:
\begin{eqnarray}\label{RF1}
f_\a(z)f^-_\a(w)&=\frac{q^{-1}z-qw}{qz-q^{-1}w}f^-_\a(w)f_\a(z)+
\frac{z(q-q^{-1})}{qz-q^{-1}w}(1+q^2)f^+_\a(q^2z)f_\a(z)\,,%F_\a(z)
\\%\ee
%\be\label{new-f}
%F_\a(z)=(1+q^2)f^+_\a(q^2z)f_\a(z)
%\ee
%and
%\be
\label{RF111}
f_{\a+\b}(z)f^-_\a(w)&=\frac{z-w}{qz-q^{-1}w}f^-_\a(w)f_{\a+\b}(z)+
\frac{(q-q^{-1})z}{qz-q^{-1}w}f_{\a+\b}(z)f_\a^-(z)\, .
\end{eqnarray}%\ee
\end{lem}
\noindent
{\it Proof.}\ \  These relations follow from applying
the integral transformation $-\oint \frac{du}{u}\frac{1}{1-w/u}$ to
the relations
\be\label{RF0}
f_\a(z)f_\a(u)=\frac{q^{-1}z-qu}{qz-q^{-1}u}f_\a(u)f_\a(z)+
(q^{-2}-q^2)\delta(q^2z/u) f_\a(q^2z)f_\a(z)\ ,
\ee
$$
f_\a(u)f_{\a+\b}(z)=\frac{q^{-1}u-qz}{u-z}f_{\a+\b}(z)f_\a(u)\,.
  $$
\hfill{$\Box$}

{}To prove Statement (ii), we use relation \r{pl-r2} and the normally
ordered expansion
$$ f_{\a+\b}(z_k) =P\sk{f_{\a+\b}(z_k})+
(q^{-1}-q)\sk{f^-_{\b}(q^{-1}z_k)f_{\a}(z_k)}^+ - f^-_{\a+\b}(z_k)\,
$$
substituted in the left-hand side of relation \r{RF12}. We then
inductively move the currents $f^-_{\b}(q^{-1}z_k)$ and
$f^-_{\a+\b}(z_k)$ to the left, using relations \rf{ffpm2} and
observing that $P(f^-_{\b}(z_k)\cdot F)=P(f^-_{\a+\b}(z_k)\cdot
F)=0$ for any element $F\in U_F$. \hfill{$\Box$}

\medskip
 We now use Statement (ii) in Proposition \ref{hyper}. It states
that the string $\ f_{\a+\b}(t_1) \cdots  f_{\a+\b}(t_k)$
$ f_{\a}(t_{k+1}) \cdots  f_{\a}(t_a)$ has simple zeroes at hyperplanes
$t_a=t_i$, $i=1,\ldots,a-1$. Substituting these conditions in \rf{RF11} and \rf{RF12}
 gives a systems of $a-1$ linear equations over the field of rational
functions $\CC(t_1,...,t_{a-1})$ for the operators $X_j(t_1,\ldots,t_{a-1})$:
\be\label{system}
\sum_{j=1}^{a-1}\frac{X_j(t_1,\ldots,t_{a-1})} {t_i-q^2t_j}=
X\cdot f^+_\a(t_i)\,, \quad i=1,..., a-1\, ,\ee
where $X={P}\sk{f_{\a+\b}(t_1)\cdots f_{\a+\b}(t_k)
f_{\a}(t_{k+1})\cdots f_{\a}(t_{a-1})}$.
The determinant of the matrix $B_{i,j}=(t_i-q^2t_j)^{-1}$ of this system is nonzero in
 $\CC(t_1,...,t_{a-1})$,
$$\det(B)=(-q^2)^{\frac{a(a-1)}{2}}\frac{\prod_{i\not=j}(t_i-t_j)^2}
{\prod_{i,j}(t_i-q^2t_j)}\,,$$
and the system hence has a unique solution over  $\CC(t_1,...,t_{a-1})$.
This implies that the operators $X_j$ are linear combinations over
$\CC(t_1,...,t_{a-1})$ of the operators $X\cdot f_\a^+(t_j)$, $j=1,...,a-1$, and
the projection of the string can therefore be represented as
\be\label{4.21}
%\begin{array}{l}
\ds {P}\sk{f_{\a+\b}(t_1)\cdots f_{\a+\b}(t_k)
f_{\a}(t_{k+1})\cdots f_{\a}(t_a)}=%\\[3mm]
%\ds\quad  = Y
X\cdot f^+_\a(t_a)-
\sum_{j=1}^{a-1}\varphi_{t_j}(t_a;t_1,\ldots,t_{a-1})X\cdot
f^+_\a(t_j)\,,
%\end{array}
\ee
%\frac{A_j(t_a;t_1,\ldots,t_{a-1})}
%{\prod_{m=1}^{a-1}(t_a-q^2t_m)}
%\Pfp\sk{f_{\a+\b}(t_1)\cdots f_{\a+\b}(t_k)
%f_{\a}(t_{k+1})\cdots f_{\a}(t_{a-1})}\Pfp\sk{f_\a(t_j)}
where
$\varphi_{t_j}(t_a;t_1,\ldots,t_{a-1})$ $=$
${A_j(t_a;t_1,\ldots,t_{a-1})}/
{\prod_{m=1}^{a-1}(t_a-q^2t_m)}$ are rational functions  whose
nomerators are polynomials in $t_a$ of degree less then $a-1$.
System \rf{system} is satisfied if the rational functions
$\varphi_{t_j}(t_a;t_1,\ldots,t_{a-1})$ have the property

$$\varphi_{t_j}(t_i;t_1,\ldots,t_{a-1})=\delta_{i,j}\,,\qquad i,j=1,...,a-1\, .$$
This interpolation problem has a unique solution given by formula \r{RF15}.

Relation \rf{4.21} now appears as a recursive relation between
projections of strings of different lengths.
The corresponding relation for the strings $f_{\a+\b}(t_1)\cdots f_{\a+\b}(t_{k})$
looks the same,
$$
\ds {P}\sk{f_{\a+\b}(t_1)\cdots f_{\a+\b}(t_a)}=
X'\cdot P(f_{\a+\b})(t_a)+
\sum_{j=1}^{a-1}\varphi_{t_j}(t_a;t_1,\ldots,t_{a-1})
X'\cdot P(f_{\a+\b})(t_j)\,,
$$
where $X'= {P}\sk{f_{\a+\b}(t_1)\cdots f_{\a+\b}(t_{a-1})}$
and the rational functions $\varphi_{t_j}(t_a;t_1,\ldots,t_{a-1})$ are
 given by relation \rf{RF15}.
  Successively applying the recursive relations and relation \rf{oop}
 gives the statement of Theorem~\ref{prstring}.
\hfill{$\Box$}

\setcounter{equation}{0}
\section{Current adjoint action and symmetrization}\label{adjoint-act}

\subsection{Projections and analytic continuation}\label{razdel5.1}

We recall \cite{DKhP} that in the completed algebra $\UqtriD$, any
product of currents $f_{\i(1)}(t_1)\cdots\!\! f_{\i(n)}(t_n)$ can be
considered an analytic function in a domain
$|t_1|\gg|t_2|\gg\cdots\gg|t_n|$, admitting an analytic continuation
to a meromorphic function in $(\CC^*)^n$. Because of commutation
relations \rf{1}, the analytic continuation of this product to the
domain
 $|t_{\si(1)}|\gg|t_{\si(2)}|\gg\cdots\gg|t_{\si(n)}|$ for any
$\sigma\in S_{n}$ is given by the series
\bn\label{ancon}
\prod\limits_{k>l\atop \si(k)<\si(l)}\frac{q^{(\i(\si(k)),\ \i(\si(l)))} -
\frac{\ds t_{\si(l)}}{\ds t_{\si(k)}}}
{1-q^{(\i(\si(k)),\ \i(\si(l)))}
\frac{\ds t_{\si(l)}}{\ds t_{\si(k)}}}
f_{\i(\si(1))}(t_{\si(1)})\cdots f_{\i(\si(n))}(t_{\si(n)}) \, .
\ee

Let $\ {}^\sigma P\sk{f_{\i(1)}(t_1)\cdots f_{\i(n)}(t_n)}\ $ denote the analytic
continuation of the projection of the current product
$P\sk{f_{\i(\si(1))}(t_{\si(1)})\cdots f_{\i(\si(n))}(t_{\si(n)})}$ from
the domain $|t_{\si(1)}|\gg|t_{\si(2)}|\gg\cdots\gg|t_{\si(n)}|$ to the domain
 $|t_1|\gg|t_2|\gg\cdots\gg|t_n|$.
\begin{prop} ${}$ \label{anprop}

\noindent {(i)} The projections $P^\pm$ commute with the analytic continuation.

\noindent {(ii)}
We have the identity of formal series in
$U\{t_{1},...,t_{n}\}$ (see \rf{Wser})
\bn
{}^\sigma P^\pm\sk{f_{\i(1)}(t_1)\cdots f_{\i(n)}(t_n)}=
\prod\limits_{k<l \atop \si^{-1}(k)>\si^{-1}(l)}
\frac{q^{(\i(k),\ \i(l))} -
t_l/t_k}%%\frac{\ds t_{l}}{\ds t_{k}}}
{1-q^{(\i(k),\ \i(l))}
t_l/t_k}%\frac{\ds t_{l}}{\ds t_{k}}}
 P^\pm\sk{f_{\i(1)}(t_1)\cdots f_{\i(n)}(t_n)}\,.
\ee
\end{prop}
{\it Proof.} Statement (i) is based on the fact that the projection operator
preserve the normal ordering in the algebra $U_F$
 (with respect to the action in the category
% $\mathcal{O}$
of highest-weight representations); in other words, it is continuous in the topology
that  defines the completion  $\overline{U}_F$. Statement (ii) follows from (i) and \rf{ancon}.
\hfill$\Box$

This proposition provides a powerful tool for the computing weight
functions. The crucial point in its application is that in contrast
to a product of currents, a projection of a current product admits a
simple analytic continuation, equivalent to the analytic
continuation of rational functions.

{\it Example.} We consider $P\sk{f_\a(t_1)f_\a(t_2)}$.
%% By proposition \ref{anprop}, we have an equallity
We have (see Sec.~\ref{sect-ex})
\begin{eqnarray}\label{ex1}
\ds{P}\sk{f_\a(t_1)f_{\a}(t_2)}&=&\ds
f_{\a}^+(t_1)f^+_\a(t_2)-\frac{(q-q^{-1})t_1}{qt_1-q^{-1}t_2}
\sk{f^+_\a(t_1)}^2\,, \\ \label{ex2}
\ds{P}\sk{f_\a(t_2)f_{\a}(t_1)}&=&\ds
f_{\a}^+(t_2)f^+_\a(t_1)-\frac{(q-q^{-1})t_2}{qt_2-q^{-1}t_1}
\sk{f^+_\a(t_2)}^2\,.
\end{eqnarray}
Equality \rf{ex1} is an equality of formal series in the domain
$|t_1|\gg |t_2|$, which means that the rational function
$t_1/(qt_1-q^{-1}t_2)$ is expanded in a power series in $t_2/t_1$;
equality \rf{ex2} is an equality of formal series
in the domain
$|t_2|\gg |t_1|$, which means that the rational function
$t_2/(qt_2-q^{-1}t_1)$ is expanded in a power series in $t_1/t_2$.
The analytic continuation to the domain $|t_1|\gg |t_2|$
in the right-hand side of \rf{ex2} amounts to the analytic continuation
of the rational function $t_2/(qt_2-q^{-1}t_1)$, which should now be expanded
in a power series in $t_2/t_1$. Therefore, the equality
$$
{}^{(12)}{P}\sk{f_\a(t_1)f_{\a}(t_2)}=
\frac{q^{2}-t_2/t_1}{1-q^{2}t_2/t_1}{P}\sk{f_\a(t_1)f_{\a}(t_2)}$$
implies a relation between formal series in the domain $|t_2|\gg
|t_1|$, \be \label{hcc1}
f_{\a}^+(t_2)f^+_\a(t_1)-\frac{(q-q^{-1})t_2}{qt_2-q^{-1}t_1}
\sk{f^+_\a(t_2)}^2 =\frac{q^{}t_1-q^{-1}t_2}{q^{-1}t_1-q^{}t_2}
f_{\a}^+(t_1)f^+_\a(t_2)-\frac{(q-q^{-1})t_1}{q^{-1}t_1-q^{}t_2}
\sk{f^+_\a(t_1)}^2  \ . \ee This is one of the basic relations in
the Borel subalgebra of $\Uqdva$. It also holds in the domain
$|t_1|\gg |t_2|$ and can be generalized to a multiple product (see
the corollary~\ref{inverse} to Theorem~\ref{sl2prop}).

\medskip

According to our definition of $q$-symmetrization (see Sec. 3.2),
the statements in Proposition~\ref{anprop}, and the rule of analytic
continuation of current product \rf{ancon},
 projections of current products with the same
 simple-root index, as well as the current products themselves, are $q$-symmetric:
\be\label{q-sym-pr}
\begin{array}{rcl}
\ds f_\a(t_1)\cdots f_\a(t_n)&=&\frac{1}{n!}\
{\overline{\rm Sym}}^n_t f_\a(t_1)\cdots f_\a(t_n)\, ,\\
\ds P^\pm\sk{f_\a(t_1)\cdots f_\a(t_n)}&=&\frac{1}{n!}\
{\overline{\rm Sym}}^n_t P^\pm\sk{f_\a(t_1)\cdots f_\a(t_n)}\, .
\end{array}\ee
We now use these arguments to write a symmetrized version
of canonical decomposition \rf{addF} of a  product of currents.

\begin{prop}\label{decprop} There is an equality of formal  series
 in the domain $|t_1|\gg |t_2|\gg ... \gg |s_{b}|$:
\bn
\begin{array}{c}
f_\a(t_1)\cdots f_\a(t_{a})f_\b(s_1)\cdots f_\b(s_{b})=
\qquad\qquad\qquad\\
=\sum_{0\leq m\leq a\atop 0\leq k\leq b}
\frac{1}{m!(a-m)!k!(b-k)!}
\ \overline{\rm Sym}^{a}_{t}\
\overline{\rm Sym}^{b}_{s}\ \Big(
\prod_{m+1\leq\ell\leq a\atop 1\leq\ell'\leq k}
\frac{qt_\ell-s_{\ell'}}{t_\ell-qs_{\ell'}}\times\\
\times \ P^-\sk{f_\a(t_1)\cdots\! f_\a(t_{m})f_\b(s_1)\cdots\!
f_\b(s_{k})}\cdot
P^+\sk{f_\a(t_{m+1})\cdots\! f_\a(t_{a})f_\b(s_{k+1})\cdots\!
f_\b(s_{b})}\Big).
\end{array}
\ee
\end{prop}
In particular, for currents of the same type, we have
\be\label{dec-ff11}
f_\b(s_{1})\cdots f_\b(s_{b})=
 \sum_{0\leq k\leq b}
\frac{1}{k!(b-k)!}\
\overline{\rm Sym}^{b}_{s}\left(
\Pfm\sk{
 f_\b(s_{1})\cdots f_\b(s_{k})     }
\cdot
\Pfp\sk{
f_\b(s_{k+1})\cdots
f_\b(t_{b})}\right).
\ee

\noindent {\it Proof}.  Directly applying expansion
\rf{addF} for the product of currents for formal series in the domain
$|t_1|\gg\cdots\gg|t_n|$ gives the relation
\be\label{dec-ff}
\begin{array}{c}
 f_{\i(1)}(t_1)\cdots\! f_{\i(n)}(t_n)=\sum_{J\subset I}
\prod_{\ell<\ell'\atop \ell\in J,\ \ell'\not\in J}
\frac{t_\ell-q^{(\i(\ell),\i(\ell'))}t_{\ell'}}
{q^{(\i(\ell),\i(\ell'))}t_\ell-t_{\ell'}}\times \\
\times\ P^-\sk{f_{\i(j'_1)}(t_{j'_1})\cdots f_{\i(j'_k)}(t_{j'_k}) }\cdot
P^+\sk{f_{\i(j''_1)}(t_{j'_1})\cdots f_{\i(j'_k)}(t_{j''_l}) }\, ,
\end{array}\ee
where the set $I=\{1,...,n\}$, its subset $J=\{j'_1,...,j'_k\}$, where $j'_1<\cdots<j'_k$,
and $I\setminus J=\{j''_1,...,j''_l\}$, where $j''_1<\cdots<j''_l$.
Applying the symmetrization procedure based on Proposition~\ref{anprop}
gives (5.7).
 \hfill$\Box$

\subsection{Current adjoint action}
\label{ad-cal}

Serre relations \rf{serre2} and \rf{serre1} admit different
representations. In Sec.~\ref{compcurrent},
 we represented them  as properties of composite currents
and strings. Here, we reformulate the Serre relations via the
current adjoint action (see \rf{adcurrent}).

\begin{lem}\label{Serre-new} Serre relations \rf{serre1} for $i=\b$ and $j=\a$ can be written as
\be\label{ad-ser1}
{\rm ad}^{(D)}_{f_\b(s_1)f_\b(s_2)}\sk{f_\a(t)}=0\,. %\qquad
%{\rm ad}^{(D)}_{f_\a(t_1)f_\a(t_2)}\sk{f_\b(s)}=0
\ee
\end{lem}
\noindent
{\it Proof}.\ \ The statement
 follows from the chain of equalities
$$
\begin{array}{rcl}
\ds {\rm ad}^{(D)}_{f_\b(s_1)f_\b(s_2)}\sk{f_\a(t)}\!&=&\ds
{\rm ad}^{(D)}_{f_\b(s_2)}\sk{
{\rm ad}^{(D)}_{f_\b(s_1)}\sk{f_\a(t)}}=
\ds {\rm ad}^{(D)}_{f_\b(s_2)}\sk{f_{\a+\b}(t)}\delta(t/qs_1)=\\[3mm]
&=&\ds \sk{f_{\a+\b}(t)f_{\b}(s_2)-f_{\b}(s_2)\psi_\b^+(s_2)^{-1}f_{\a+\b}(t)
\psi_\b^+(s_2)}\delta(t/qs_1)=\\[3mm]
&=&\ds \sk{f_{\a+\b}(t)f_{\b}(s_2)-f_{\b}(s_2)f_{\a+\b}(t)
\frac{t-q^3s_2}{qt-q^2s_2}}\delta(t/qs_1)=0\,.
\end{array}
$$
%based on the commutation relations \r{com-f}, \r{com-bb} and
%\r{3}.
\hfill$\Box$

 This lemma implies the following proposition.
\begin{prop}\label{prop5.4}
For $a\geq k$, the identity of formal series
\be\label{sym-s1}
\begin{array}{l}
\ds P\sk{{\rm ad}^{(D)}_{f_\b(u_1)\cdots f_\b(u_k)}\sk{f_\a(t_1)\cdots
f_\a(t_a)}}= \frac{1}{k!(a-k)!}\\
\qquad\ds \times\  \overline{\rm Sym}^{a}_{t}\ \left(
 P\sk{f_\a(t_1)\cdots f_\a(t_{a-k}) f_{\a+\b}(t_{a-k+1})\cdots
 f_{\a+\b}(t_{a})}\phantom{\prod_{i=1}^b}
 \right.\\
\left.\qquad\ds\times \prod_{i<j}^{k}
\frac{q^{-1}\tilde{t}_{i}-q^{}\tilde{t}_j}
{\tilde{t}_i-\tilde{t}_j}\ \widetilde{\rm Sym}^{k}_{\tilde{t}}
\sk{\prod_{i=1}^k \delta\sk{\frac{\tilde{t}_i}{qu_i}}  }  \right)
\end{array}
\ee
holds in the domain $|t_1|\gg\cdots\gg|t_a|$, where $\tilde{t}_i=t_{a-k+i}$, $i=1,\ldots,k$.
\end{prop}
\noindent
{\it Proof.}\ \ The proof is a  combinatorial exercise involving with definition of
current adjoint action \rf{adcurrent}, its  properties including
\rf{ad-ser1}, and  the relation
$$
{\rm ad}^{(D)}_{f_\b(s)} (F_1\cdot F_2)= F_1\cdot{\rm ad}^{(D)}_{f_\b(s)}
  (F_2)+ {\rm ad}^{(D)}_{f_\b(s)} (F_1)\cdot \psi^+_\b(s)^{-1}\, F_2\,
  \psi^+_\b(s)\,.
$$
We also note the rol of Proposition~\ref{anprop}, which allows using
the commutation relation between total currents under the projections without paying
attention to  $\d$-function terms. With this taken into account and after necessary
combinatorial rearrangments, we obtain
$$%\be\label{sym-s1-pr}
\begin{array}{l}
\ds P\sk{{\rm ad}^{(D)}_{f_\b(u_1)\cdots f_\b(u_k)}\sk{f_\a(t_1)\cdots
f_\a(t_a)}}= \frac{1}{k!(a-k)!}
\ \overline{\rm Sym}^{a}_{t}\ \overline{\rm Sym}^{k}_{u}\ \left(
\prod_{i=1}^k\prod_{\ell=1}^{i-1} \frac{qt_{a-k+i}-u_{\ell}}
{t_{a-k+i}-qu_{\ell}}\times\right.\\
\ds\left.\phantom{\prod_{i=1}^b}\quad
\times P\sk{ f_\a(t_1)\cdots f_\a(t_{a-k})
{\rm ad}^{(D)}_{f_\b(u_{1})}\sk{f_\a(t_{a-k+1})}\cdots
 {\rm ad}^{(D)}_{f_\b(u_{k})}\sk{f_\a(t_{a})}}
\right).
\end{array}
$$%\ee
Because the adjoint action in the last formula produces a product
of $\delta$-functions $\prod_i \delta(t_{a-k+i}/q u_{i})$, $i=1,\ldots,k$ (see \r{adac4}),
we can now move the rational function from under the
$q$-symmetrization with respect to the variables $u_i$ and replace the symmetrization
with the $q^{-1}$-symmetrization with respect to the variables $t_{a-k+i}=\tilde{t}_{i}$.
Proposition~\ref{prop5.4} is thus proved.
\hfill$\Box$

\subsection{Proof of Theorem~\ref{act-cal}}
Our goal is to reduce an expression
$P\sk{f_\a(t_1)\cdots f_\a(t_{a})f_\b(s_1)\cdots
f_\b(s_{b})}$
to projections of strings. In this expression, we substitute
decomposition  \r{dec-ff11} for the current product $f_\b(s_1)\cdots
f_\b(s_{b})$:
\be\label{po3}
\begin{array}{l}
\ds P\sk{f_{\a}(t_1)\cdots f_{\a}(t_{\a})f_{\b}(s_1)\cdots
f_{\b}(s_{b})}= \sum_{k=0}^{b}\ \
\frac{1}{k!(b-k)!}\times\\
\ds\quad \times \overline{{\rm Sym}}^{b}_{s}\left(
\Pfp\sk{f_{\a}(t_{1})\cdots  f_{\a}(t_{a}) \Pfm
\sk{f_{\b}(s_1)\cdots f_{\b}(s_{k})}}
\cdot \Pfp\sk{f_{\b}(s_{k+1})\cdots
f_{\b}(s_{b})}
\right).
\end{array}
\ee
 We now use a strengthened coideal property of  the
subalgebra $U_f^-$.% given in the Proposition~\ref{subalgebras}
\begin{prop}${}$

\noindent (i) For any element
$F\in U_f^-$, we have
\be\label{coiF1}
\Delta^{(D)}\, F = 1\ot F + F'\ot F'', \quad\mbox{such that}\quad
F'\in U^-_f\quad\mbox{and}\quad\varepsilon (F')=0\,.
\ee
\noindent (ii) For any product $f_{\i(1)}(t_{1})\cdots  f_{\i(n)}(t_{n})$,
we have the equality of series in $U\{t_1,...,t_n\}$ (see \r{Wser})
\be\label{po5a}
\begin{array}{l}
\ds P\sk{f_{\i(1)}(t_{1})\cdots  f_{\i(m)}(t_{m}) \Pfm
\sk{f_{\i(m+1)}(t_{m+1})\cdots f_{\i(n)}(t_{n})}}=\\ \qquad
\ds= P \sk{
{\rm ad}^{(D)}_{\Pfm
\sk{f_{\i(m+1)}(t_{m+1})\cdots f_{\i(n)}(t_{n})}}
f_{\i(1)}(t_{1})\cdots  f_{\i(m)}(t_{m}) }.
\end{array} %=\\
\ee
\end{prop}
{\it Proof.}\ \ It suffices to verify statement (i) for generators of the algebra
$U_f^-$, where it is a direct observation. Statement (ii) is a direct
consequence of (i).
\hfill$\Box$

We use a particular case of \rf{po5a},
\be\label{po5}
%\begin{array}{l}
\ds P\sk{f_{\a}(t_{1})\cdots  f_{\a}(t_{a}) \Pfm
\sk{f_{\b}(s_1)\cdots f_{\b}(s_{k})}}%=\\ \qquad
\ds= P \sk{
{\rm ad}^{(D)}_{\Pfm
\sk{f_{\b}(s_1)\cdots f_{\b}(s_{k})}}
f_{\a}(t_{1})\cdots  f_{\a}(t_{a}) }\,, %=\\
\ee
and substitute an integral representation of the projection
$P^-
\sk{f_{\b}(s_1)\cdots f_{\b}(s_{k})}$ in it:
\be\label{intt2}
\!P^-\sk{f_\b(s_1)\cdots f_\b(s_k)}=
\prod_{i<j}\frac{s_i-s_j}{qs_i-q^{-1}s_j}
\mathop{\oint\!\cdots\!\oint}%\limits_{ w_a\gg \cdots \gg w_1}
Y(^\op\overline{u};\,^\op\overline{s})f_\b(u_1)\frac{du_1}{u_1}\cdots f_\b(u_k)
\frac{du_k}{u_k}\,.
\ed
Then the right-hand side of equality \r{po5} becomes
\begin{equation}\label{5.16}
\begin{array}{l}
\ds
\prod_{i<j}^k \frac{s_i-s_j}{qs_i-q^{-1}s_j}
\oint\!\cdots\!\oint \prod_{i=1}^k \frac{du_i}{u_i} \
Y(^\op\overline{u},\,^\op\overline{s}) \
P\sk{{\rm ad}^{(D)}_{f_\b(u_1)\cdots f_\b(u_k)}\sk{f_\a(t_1)\cdots
f_\a(t_a)}}=\\
\ds\quad=\ \frac{1}{k!(a-k)!}\ \ \overline{\rm Sym}^a_t\left(
P\sk{f_\a(t_1)\cdots f_\a(t_{a-k})f_{\a+\b}(t_{a-k+1})\cdots f_{\a+\b}(t_a)}
\phantom{\prod_{i<j}}\right.\\
\ds\quad \times\left.
\prod_{i<j}\frac{s_i-s_j}{qs_i-q^{-1}s_j}
\frac{q^{-1}\tilde{t}_{i}-q^{}\tilde{t}_j}
{\tilde{t}_i-\tilde{t}_j}\
\widetilde{\rm Sym}^k_{\tilde{t}}\sk{
{\oint\!\cdots\!\oint}
\prod_{i=1}^k \frac{du_i}{u_i}\,\delta \sk{\frac{\tilde{t}_i}{qu_{i}}}
Y(^\op\overline{u};\,^\op\overline{s})}\right)\, ,
\end{array}
\end{equation}
where $\tilde{t}_i=t_{a-k+i}$, $i=1,\ldots,k$.
After the integration, the last line in
  \rf{5.16} becomes
\begin{equation}\label{5.17}
\begin{array}{l}
\ds \prod_{i<j}\frac{s_i-s_j}{qs_i-q^{-1}s_j}
\frac{q^{-1}\tilde{t}_{i}-q^{}\tilde{t}_j}
{\tilde{t}_i-\tilde{t}_j}\
\widetilde{\rm Sym}^k_{\tilde{t}}\sk{
Y(q^{-1}\cdot({}^\op\tilde{t});\,^\op\overline{s})}\ = \\
\ds =\ \prod_{i<j}\frac{s_i-s_j}{qs_i-q^{-1}s_j}
\frac{q^{}\tilde{t}_{i}-q^{-1}\tilde{t}_j}
{\tilde{t}_i-\tilde{t}_j}\
\overline{\rm Sym}^k_{\tilde{t}}\sk{
Y(q^{-1}\cdot\tilde{t};\,^\op\overline{s})}\ = \
\overline{\rm Sym}^k_{s}\sk{
Y(q^{-1}\cdot\tilde{t};\,\overline{s})}\ .
\end{array}
\end{equation}
In the first equality in  \rf{5.17}, relation \rf{qs3r} between
the $q$- and $q^{-1}$-symmetrizations
is used,  while the second equality is obtained from
combinatorial identity   \rf{idZZ1}.

Projection \r{po5} is now given by
\be\label{xx4}
\begin{array}{l}
\ds P\sk{
{\rm ad}^{(D)}_{\Pfm
\sk{f_{\b}(s_1)\cdots f_{\b}(s_{k})}}
f_{\a}(t_{1})\cdots  f_{\a}(t_{a}) }=\frac{1}{k!(a-k)!}\\
\ds\qquad\times\ \overline{\rm Sym}_{t}^{a}\ \overline{\rm Sym}_{s}^{k}\
\Big(
 Y(q^{-1}t_{a-k+1},\ldots,q^{-1}t_{a};s_1,\ldots,s_k)\\
 \ds\qquad\qquad\times P\sk{
f_{\a}(t_{1})\cdots
f_{\a}(t_{a-k})
f_{\a+\b}(t_{a-k+1})\cdots f_{\a+\b}(t_{a})}\Big).
\end{array}
\ee

We now return to the proof of Theorem~\ref{act-cal}. By definition
of $q$-symmetrization \r{qs1ra}, the right-hand side of the last
formula is $q$-symmetric with respect to the variables
$s_1,\ldots,s_{k}$. Additional $q$-symmetrization then cancels the
unwanted factorial (see \r{q-sym-pr}), and we obtain the proof of
Theorem~\ref{act-cal}. \hfill$\Box$

%\section{Discussion}

\section*{Acknowledgement}

Many ideas of this paper appeared during one
author's (S.P.) stay at the Mathematical
Department of Kyushu University in the autumn 2003, and an
essential part of this work was
done during both authors' visit to the Max Plank Institut f\"ur Mathematik in the spring
2004. The authors thank these scientific centers for a stimulating
scientific atmosphere.

This paper was supported in part by the INTAS (Ggrant No.
OPEN-03-51-3350), the Heisen\-berg-Landau program, the Russian
Foundation for Basic Research (Grant No. 04-01-00642), and the
Program for Supporting Leading Scientific Schools  (Grant No.
NSh-1999.2003.2.

\app{The projection operator $P^-$}

\noindent
In this appendix, we collect the most important  formulas for the opposite projection
operator $\Pfm$.

{\bf 1.} Just as for the projection $P=P^+$, an expression
$P^-\left(f_{\i(i_1)}(t_{i_1})\cdots f_{\i(i_n)}(t_{i_n})
\right)$ is a series in the domain  $|t_{i_1}|\gg\cdots\gg|t_{i_n}|$
 admitting an analytic continuation to different asymptotic zones.
In the sense of  analytic continuation, it has the same properties
as the analogous expressions for the projections  $P^+$.

We recall the formula for the projection $\Pfp$ of the product of currents
corresponding to the same root:
$$\Pfp\sk{f_\a(t_1)\cdots f_\a(t_a)}=
f^+_\a(t_1)f^+_\a(t_2;t_1)\cdots
f^+_\a(t_{a-1};t_1,\ldots,t_{a-2})f^+_\a(t_a;t_1,\ldots,t_{a-1}).$$
The structure of this formula was explained in the proof of the Theorem~\ref{prstring}.
From analogous considerations (see the next appendix), we can find that a
similar triangular decomposition holds for the projection
 $\Pfm$ of the product of currents with the same root index,
\be\label{dec-m} \Pfm\sk{f_\b(s_1)\cdots f_\b(s_b)}=(-1)^b
f^-_\b(s_1;s_2,\ldots,s_{b})f^-_\b(s_{2};s_3,\ldots,s_{b}) \cdots
f^-_\b(s_{b-1};s_b)f^-_\b(s_b), \ee where the currents
$f^-_\b(s_{k};s_{k+1},\ldots,s_{b})$ are defined as sums
\be\label{long-cur-m} f^-_\b(s_k;s_{k+1},\ldots,s_b)=f^-_\b(s_k)-
\sum_{m=k+1}^b
\tilde{\varphi}_{s_m}(s_k;s_{k+1},\ldots,s_b)f^-_\b(s_m) \ee and the
rational functions, $\tilde{\varphi}_{s_j}(s;s_1,\ldots,s_{b})$
\be\label{RF15m}
\tilde{\varphi}_{s_j}(s;s_1,\ldots,s_{b})=\prod_{i=1,\ i\neq j}^{b}
\frac{s-s_i}{s_j-s_i}\prod_{i=1}^{b}\frac{qs_j-q^{-1}s_i}{qs-q^{-1}s_i}\,
, \ee as functions of the variable $s$ have simple poles at the
points $s=q^{-2}s_i$, $i=1,\ldots,b$, tend to zero as  $s\to\infty$,
and have the property
$\varphi_{s_j}(s_i;s_1,\ldots,s_b)=\delta_{ij}$. These conditions
define the set of functions (\ref{RF15m}) uniquely.

As can the projection  $\Pfp$, projection \r{dec-m} for the product of currents can be
written in the reverse order
$$
\Pfm\sk{f_\b(s_1)\cdots f_\b(s_b)}=(-1)^b \prod_{1\leq i<j\leq b}
\frac{q^{-1}s_i-qs_j}{qs_i-q^{-1}s_j}\
f^-_\b(s_b;s_{b-1},\ldots,s_{1})
\cdots f^-_\b(s_{2};s_1)f^-_\b(s_1)\ .
$$

Expression  \r{dec-m} can also be written as an integral transform
of the product of the total currents,
$$
P^-\sk{f_\b(s_1)\cdots f_\b(s_b)}=
\prod_{i<j}\frac{s_i-s_j}{qs_i-q^{-1}s_j}
\oint
\prod_{k=1}^b \frac{du_k}{u_k}Y( ^\op\overline{u},^\op \overline{s})
f_\b(u_1)\cdots f_\b(u_b)= $$
$$=\prod_{i<j}\frac{s_i-s_j}{qs_i-q^{-1}s_j}
\oint
\prod_{k=1}^b \frac{du_k}{u_k} \frac{1}{1-s_k/u_k}\prod_{i=k+1}^b
\frac{q-q^{-1}s_i/u_k}{1-s_i/u_k}
f_\b(u_1)\cdots f_\b(u_b)\ ,
$$
which was already used in the preceding section in the proof of Theorem~\ref{act-cal}.

{\bf 2.}
Calculation the projection  $P^-(f_\a(t_1)\cdots f_\a(t_{a})f_\b(s_1)\cdots
f_\b(s_{b}))$ also reduces to calculating the projections of strings,
\be\label{po2}
\begin{array}{l}
\ds P^-\sk{f_\a(t_1)\cdots f_\a(t_{a})f_\b(s_1)\cdots
f_\b(s_{b})}=\\
\ds\  = \!\!\sum_{k=0}^{\min\{a,b\}}\!\!\!\!
\frac{1}{k!(a-k)!(b-k)!}
\ \overline{{\rm Sym}}^{a}_{t}\overline{{\rm Sym}}^{b}_{s}\left(
P^-\sk{f_\a(t_{1})\cdots  f_\a(t_{a-k})}
\right. \times\\
\ds\qquad  \!\! P^-\sk{
f_{\a+\b}(qs_{1})\cdots f_{\a+\b}(qs_{k})
f_\b(s_{k+1})\cdots
f_\b(s_{b})}
\left.Z(q^{-1}t_{a-k+1},...,q^{-1}t_{a};s_1,\ldots,s_k)
\right),
\end{array}
\ee
where the series $Z(\overline{t},\overline{s})$ is defined in \rf{rat-Y} and
the projection of a string is given by
\be\label{po44}
\begin{array}{l}
\ds \Pfm\sk{f_{\a+\b}(qs_1)\cdots f_{\a+\b}(qs_{k})f_\b(s_{k+1})
\cdots f_\b(s_{b})}=\\
\ds\quad =
\prod_{1\leq i<k+1\leq j\leq b}\frac{qs_i-q^{-1}s_j}{s_i-s_j}
\prod_{1\leq i<j\leq b}\frac{q^{-1}s_i-qs_j}{qs_i-q^{-1}s_j}\\
\ds\qquad\times\ \Pfm\sk{f_\b(s_b;s_{b-1},\ldots,s_1)}
\cdots \Pfm\sk{f_\b(s_{k+1};s_{k},\ldots,s_1)}\\
\ds\qquad\quad\times\ \Pfm\sk{f_{\a+\b}(qs_k;qs_{k-1},\ldots,qs_1)}\cdots
\Pfm\sk{f_{\a+\b}(qs_1)}\ .
\end{array}
\ee
In this formula, the currents $f_\gamma(t_i;t_{i-1},\ldots,t_1)$ for the roots  $\gamma=\b,\a+\b$
are defined by relations \r{long-cur-m}
with coefficient functions  \r{RF15m}, which are invariant under a simultaneous
scaling of all variables. Single current projections are defined by the formulas
 \rf{3.14} such that
 $P^-\sk{f_{\a}(t)}=-f^-_\a(t)$,
$ P^-\sk{f_{\a+\b}(qs)}=
-\big(f_\a[0]f^-_\b(s)-q^{-1}f^-_\b(s)f_\a[0]\big)$, and
$P^-\sk{f_{\b}(s)}=-f^-_\b(s)$.
In particular, we have
$$
\begin{array}{rcl}
\ds{P}^-\sk{f_\a(t_1)f_{\a}(t_2)}&=&\ds
\sk{f^-_\a(t_1)-\frac{(q-q^{-1})t_2}{qt_1-q^{-1}t_2}
f^-_\a(t_2)}f_{\a}^-(t_2)= \\
\ds &=&\ds\frac{q^{-1}t_1-qt_2}{qt_1-q^{-1}t_2}
\sk{f^-_\a(t_2)-\frac{(q-q^{-1})t_1}{qt_2-q^{-1}t_1}
f^-_\a(t_1)}f_{\a}^-(t_1),\\
\ds{P}^-\sk{f_{\a+\b}(qs_1)f_\b(s_2)}&=&-\ds\frac{q^{-1}s_1-qs_2}{s_1-s_2}\
\sk{f_{\b}^-(s_2)-\frac{(q-q^{-1})s_1}{qs_2-q^{-1}s_1}
f_{\b}^-(s_1)}P^-(f_{\a+\b}(qs_1))\ .
\end{array}$$

{\bf 3.} The proof of formula \r{po2} is based on one more
reformulation of the Serre relations involving the right adjoint
action $\tilde{\rm ad}_x^{(D)}(y)$ (see \r{adcurrent}),
\be\label{ad-ser2} \tilde{\rm
ad}^{(D)}_{f_\a(t_1)f_\a(t_2)}\sk{f_\b(s)}=0\, , \ee and on an
identity analogous to the one proved in Proposition~\ref{prop5.4}.
Namely, for $k\leq b$, \be\label{sym-s2}
\begin{array}{l}
\ds \Pfm\sk{\tilde{\rm ad}^{(D)}_{f_\a(t_1)\cdots f_\a(t_k)}
\sk{f_\b(s_1)\cdots
f_\b(s_b)}}= \frac{1}{k!(b-k)!}
\overline{\rm Sym}^{b}_{s}\overline{\rm Sym}^{k}_{t}
\left(\prod_{i=1}^k\prod_{\ell=i+1}^{k}
\frac{qt_{\ell}-s_i}
{t_{\ell}-qs_i}\right.\times\\
\ds \left.\phantom{\prod_{i=1}^k}\quad\times\
\Pfm\sk{\tilde {\rm ad}^{(D)}_{f_\a(t_1)}\sk{f_\b(s_{1})}\cdots
\tilde {\rm ad}^{(D)}_{f_\a(t_k)}\sk{f_\b(s_{k})}
 f_\b(s_{k+1})\cdots  f_\b(s_{b})
}\right)\,.
\end{array}
\ee

\app{Commutation relations
 with projections
of currents}

\noindent Commutation relations \r{com-f} imply the rules for
moving the half-currents $f^-(z)$ to the left
and the half-currents $f^+(z)$ to the right through the total currents.
\begin{prop}\label{h-t-com}
We have the equalities
\be\label{t-f-p}
\begin{array}{rcl}
f^+_{\a}(z)f_{\b}(w)&=&\frac{qz- w}{z-qw}
f_{\b}(w)f^+_{\a}(z)+\\
&+&\frac{qw(q^{-1}-q)}{z-qw}
f_{\b}(w)f^+_{\a}(qw)+\frac{qw}{z-qw}
f_{\a+\b}(qw)\,,
\end{array}
\ee
\be\label{t-f-m}
\begin{array}{rcl}
f_{\a}(z)f^-_{\b}(w)&=&\frac{qz-w}{z-qw}
f^-_{\b}(w)f_{\a}(z)+\\
&+&\frac{z(q^{-1}-q)}{z-qw}
f^-_{\b}(q^{-1}z)f_{\a}(z)-\frac{z}{z-qw}
f_{\a+\b}(z)\,.
\end{array}
\ee
\end{prop}
\noindent
{\it Proof}.\ The proof of relation (\ref{t-f-p}), for example,
 is based on the decomposition of the kernel
$$
\frac{qu-w}{(u-qw)(z-u)}=
\frac{qz-w}{z-qw}\ \frac{1}{z-u}+
\frac{ qw(q^{-1}-q)}{z-qw}\ \frac{1}{qw-u}
$$
into the sum of two kernels.\hfill$\Box$

Analogously, from relations \r{com-aa}, we have the following proposition.
\begin{prop}\label{h-t-S} The nontrivial commutation relations between the
  total and half-currents  are
\be\label{ffpm1}
\begin{array}{rcl}
f^+_{\a}(z)f_{\a+\b}(w)&=&\frac{q^{-1}z-qw}{z-w}f_{\a+\b}(w)f^+_{\a}(z)+
\frac{(q-q^{-1})w}{z-w}f_{\a+\b}(w)f^+_{\a}(w)\,,\\
f_{\a}(z)f^-_{\a+\b}(w)&=&\frac{q^{-1}z-qw}{z-w}f^-_{\a+\b}(w)f_{\a}(z)+
\frac{(q-q^{-1})z}{z-w}f^-_{\a+\b}(z)f_{\a}(z)
\end{array}
\ee
and
\be\label{ffpm2}
\begin{array}{rl}
f^+_{\a+\b}(qz)f_{\b}(w)&=
\frac{q^{-1}z-qw}{z-w}
f_{\b}(w)f^+_{\a+\b}(qz)+
\frac{(q-q^{-1})w}
{z-w}f_{\b}(w)f^+_{\a+\b}(q^{p-i}w)\,,\\
f_{\a+\b}(qz)f^-_{\b}(w)&=
\frac{q^{-1}z-qw}{z-w}
f^-_{\b}(w)f_{\a+\b}(qz)+
\frac{(q-q^{-1})z}
{z-w}f^-_{\b}(z)f_{\a+\b}(qz)\,,\\
f_{\a+\b}(z)f^-_{\a+\b}(w)&=\frac{q^{-1}z-qw}{qz-q^{-1}w}
f^-_{\a+\b}(w)f_{\a+\b}(z)+
\frac{z(q-q^{-1})}{qz-q^{-1}w}(1+q^2)f^+_{\a+\b}(q^2z)f_{\a+\b}(z)\,.
\end{array}
\ee
\end{prop}

{\it Remark.}
The meaning of the commutation relations in this proposition
is that we can move the half-currents $f^-_\gamma(w)$ to the left and the
half-currents $f^+_\gamma(z)$ to the right through the total currents such
that the total currents are unchanged and only shifted half-currents multiplied
by rational functions arise. In this paper, we often use these  properties of
exchange between the total and half-currents.

The following proposition describes the commutation relations between
projections of a com\-po\-si\-te-root current and the projection of simple-root currents.
\begin{prop}\label{restore-ord}
There are the equalities
$$
\begin{array}{l}
\ds \Pfp\sk{f_{\a+\b}(t_1)}f^+_\a(t_2)-
qf^+_\a(t_2)\Pfp\sk{f_{\a+\b}(t_1)}=\\
\ds\quad =
\frac{q-q^{-1}}{t_1-t_2}
(f^+_\a(t_1)-f^+_\a(t_2))\sk{t_1\Pfp\sk{f_{\a+\b}(t_1)}-
t_2\Pfp\sk{f_{\a+\b}(t_2)}}\,,\\[5mm]
\ds f^-_\b(s_1) \Pfm\sk{f_{\a+\b}(qs_2)}-
q^{-1} \Pfm\sk{f_{\a+\b}(qs_2)}f^-_\b(s_1) =\\
\ds\quad =
\frac{q-q^{-1}}{s_1-s_2}\sk{s_1\Pfm\sk{f_{\a+\b}(qs_1)}-
s_2\Pfm\sk{f_{\a+\b}(qs_2)}}
(f^-_\b(s_1)-f^-_\b(s_2))\,.
\end{array}
$$
\end{prop}
\noindent
{\it Proof.}\ The proof consists in applying the screening operator
$S_{f_\b[0]}$  to equality \r{hcc1} and $\tilde S_{f_\a[0]}$ to the analogous equality
where the projection of currents  $f^+_\a(t_i)$ is replaced with
$f^-_\b(s_i)$. the proof also uses  the Serre relations written in terms
of the projections of currents  as
$$
\begin{array}{l}
\ds q f^+_\a(t_1) \Pfp\sk{f_{\a+\b}(t_2)}
+ q f^+_\a(t_2) \Pfp\sk{f_{\a+\b}(t_1)}=\\[5mm]
\ds\quad =\ \Pfp\sk{f_{\a+\b}(t_1)} f^+_\a(t_2)
+\Pfp\sk{f_{\a+\b}(t_2)} f^+_\a(t_1)\,,\\[7mm]
\ds q f^-_\b(s_1) \Pfm\sk{f_{\a+\b}(s_2)}
+ q f^-_\b(s_2) \Pfm\sk{f_{\a+\b}(s_1)}=\\[5mm]
\ds\quad =\ \Pfm\sk{f_{\a+\b}(s_1)} f^-_\b(s_2)
+\Pfm\sk{f_{\a+\b}(s_2)} f^-_\b(s_1)\,.
\end{array}
$$

\hfill$\Box$
\bigskip

\app{Direct proof of Theorem~\ref{act-cal}}

The first step in this proof is the same as in \r{po3}.
We then observe that the first projection $\Pfp$ in \r{po3}
vanishes for $k>a$. This explains the upper summation limit in \r{po1}.
To prove this formula, we calculate the projection
$$%\be\label{po55}
\Pfp\sk{f_{\a}(t_{1})\cdots  f_{\a}(t_{a}) \Pfm
\sk{f_{\b}(s_1)\cdots f_{\b}(s_{k})}}
$$%\ee
or
$$
\Pfp\sk{f_{\a}(t_{1})\cdots  f_{\a}(t_{a})
f^-_{\b}(s_1;s_2,\ldots,s_{k})\cdots
f^-_{\b}(s_{k-1};s_{k})f^-_{\b}(s_{k})}\ ,
$$
moving the half-currents $f^-_{\b}(s_m;s_{m+1},\cdots,
s_{k})$, $m=1,\ldots,k$,
to the left using commutation relations \r{t-f-m}.
According to these commutation relations, the composite currents
$f_{\a+\b}(t_j)$ are created at the positions $j=1,\ldots,a$.
We next move these composite currents to the right using
commutation relations \r{com-aa}.
Here, we again use the statements in Proposition~\ref{anprop}.
After these commutations, we obtain the sum over all nonordered
subsets $J=\{j_1,\ldots,j_{k}\}\in\{1,\ldots,a\}$:
$$
\begin{array}{c}
\ds \prod_{i<j}^{k}\frac{s_i-s_j}{qs_i-q^{-1}s_j}
\ \sum_{J}
\prod_{m=1}^{k}\sk{\prod_{\ell_m=1}^{j_m-1}
\frac{qt_{j_m}-q^{-1}t_{\ell_m}}{t_{j_m}-t_{\ell_m}}
\prod_{\ell_m=j_m+1}^{a}
\frac{q^{-1}t_{j_m}-qt_{\ell_m}}{t_{j_m}-t_{\ell_m}}}
\times\\[3mm]
\ds
\times
\prod_{m=1}^{k}
\sk{\frac {t_{j_m}}{t_{j_m}-qs_m} \prod_{i=m+1}^{k}
\frac{qt_{j_m}-s_i}{t_{j_m}-qs_i} }
\Pfp\!\!\sk{f_{\a+\b}(t_{j_1})\cdots f_{\a+\b}(t_{j_{k}})
\mathop{\underbrace{
f_{\a}(t_{1})f_{\a}(t_{2})\cdots f_{\a}(t_{a-1})
f_{\a}(t_{a}})}
\limits_{{\rm currents\ depending\ on}\
 t_{j_1},\ldots,t_{j_{k}}\ {\rm omitted}}
}\!.
\end{array}
$$

The next step is to use the commutation relations \r{com-aa} to move
the group of the composite currents
$f_{\a+\b}(t_{j_1})\cdots f_{\a+\b}(t_{j_{k}})$
under projection from left to  right,
\be\label{xx1}
\begin{array}{c}
\ds \prod_{i<j}^{k}\frac{s_i-s_j}{qs_i-q^{-1}s_j}
\ \sum_{J}
\prod_{i<j\atop i,j\in J}
\frac{q^{-1}t_{i}-qt_{j}}{t_{i}-t_{j}}
\prod_{i<j\atop i\in J,\ j\not\in J}
\frac{q^{-1}t_{i}-qt_{j}}{qt_{i}-q^{-1}t_{j}}
\
\prod_{m=1}^{k}
\sk{\frac {t_{j_m}}{t_{j_m}-qs_m} \prod_{i=m+1}^{k}
\frac{qt_{j_m}-s_i}{t_{j_m}-qs_i} }
\times\\[3mm]
\ds
\times\ \Pfp\sk{
\mathop{\underbrace{
f_{\a}(t_{1})f_{\a}(t_{2})\cdots f_{\a}(t_{a-1})
f_{\a}(t_{a}})}
\limits_{{\rm currents\ depending\ on}\
 t_{j_1},\ldots,t_{j_{k}}\ {\rm omitted}}
f_{\a+\b}(t_{j_1})\cdots f_{\a+\b}(t_{j_{k}})}.
\end{array}
\ee We now use the fact that the summation in this formula ranges
nonordered sets $J$ of size $k$. This means that this summation can
be decomposed into two summations: first, over all different but
ordered choices $\{j_1<j_2<\cdots< j_{k}\}$ from the set
$\{1,2,\ldots,a\}$ and, second, over all permutations among fixed
$\{j_1,j_2,\ldots, j_{k}\}$. Because of commutation relations
\r{com-aa} between composite currents, this second summation can be
written as a $q^{-1}$-symmetrization with respect to a fixed subset
$\{j_1,j_2,\ldots, j_{k}\}$ of the function
$$
\prod_{m=1}^{k}
\sk{\frac {1}{t_{j_m}-qs_k} \prod_{i=m+1}^{k}
\frac{qt_{j_m}-s_i}{t_{j_m}-qs_i} }
$$
because the other ingredients in formula \r{xx1} are stable under
permutation of $\{j_1,j_2,\ldots, j_{k}\}$. Using combinatorial identity
\r{idZZ1}, which can be written as
$$%\be\label{xx22}
\begin{array}{c}
\ds\prod_{i<j}^{k}\frac{q^{-1}t_{i}-qt_{j}}{t_{i}-t_{j}}
\ \ \widetilde{\rm Sym}_{t}^{k}\sk{\
\prod_{m=1}^{k}
\sk{\frac {1}{t_m-qs_{m}} \prod_{i=m+1}^{k}
\frac{qt_m-s_{i}}{t_m-qs_{i}} }}
=\\
\ds=\prod_{i<j}^{k}\frac{qs_i-q^{-1}s_j}{s_i-s_j} \ \ \overline{\rm
Sym}_{s}^{k} \sk{\ \prod_{m=1}^{k} \sk{\frac {1}{t_{m}-qs_m}
\prod_{i=1}^{m-1} \frac{qt_{m}-s_i}{t_{m}-qs_i} }}\ ,
\end{array}
$$%\ee
we can now rewrite expression \r{xx1} as
\be\label{xx3}
\begin{array}{c}
\ds  \sum_{J\atop j_1<\cdots< j_{k}}
\prod_{i<j\atop i\in J,\ j\not\in J}
\frac{q^{-1}t_{i}-qt_{j}}{qt_{i}-q^{-1}t_{j}}
\
\overline{\rm Sym}_{s}^{k}\
\sk{\prod_{m=1}^{k}
\frac {t_{j_m}}{t_{j_m}-qs_k} \prod_{i=1}^{m-1}
\frac{qt_{j_m}-s_i}{t_{j_m}-qs_i} }
\times\\[3mm]
\ds
\times\Pfp\sk{
\mathop{\underbrace{
f_{\a}(t_{1})f_{\a}(t_{2})\cdots f_{\a}(t_{a-1})
f_{\a}(t_{a}})}
\limits_{{\rm currents\ depending\ on}\
 t_{j_1},\ldots,t_{j_{k}}\ {\rm omitted}}
f_{\a+\b}(t_{j_1})\cdots f_{\a+\b}(t_{j_{k}})}.
\end{array}
\ee The summation over all ordered sets $J$ in \r{xx3} can now be
written as a $q$-symmetrization with respect to all the variables
 $t_1,\cdots,t_{a}$, and we obtain \r{xx4}. Repeating the argument
 given at the end of Sec.~\ref{ad-cal}, we finish the direct proof of
Theorem~\ref{act-cal}.


\begin{thebibliography}{9999}

\bibitem{KR83} P. Kulish, N. Reshetikhin. Diagonalization of $GL(N)$
  invariant transfer matrices and quantum $N$-wave system (Lee model)
{\sl J.Phys. A: Math. Gen.} {\bf 16} (1983), L591-L596.


\bibitem{D88} V. Drinfeld. New realization of Yangians and quantum
  affine algebras. {\it Sov. Math. Dokl.} {\bf 36} (1988), 212--216.


\bibitem{E} B. Enriquez, On correlation functions of Drinfeld
currents and shuffle algebras,  {\it Transformation Groups}
{\bf 5:2} (2000), 111-120.

\bibitem{EKPT} B. Enriquez, S. Khoroshkin, and S. Pakuliak,
Weight functions and Drinfeld currents. Preprint ITEP-TH-40/50,
 math.QA/0610398.


\bibitem{ER} B. Enriquez, V. Rubtsov. Quasi-Hopf algebras associated with
$\mathfrak{sl}_2$ and complex curves. Israel J Math {\bf 112} (1999) 61--108.



\bibitem{DKh} J. Ding, S. Khoroshkin. Weyl group extention of quantized
current algebras. {\it Transformation Groups} {\bf 5} (2000) 35--59.



\bibitem{DKhP1} Ding, J.; Khoroshkin, S.; Pakuliak, S.
 Integral presentations for the universal $R$-matrix. Lett. Math.
Phys. \textbf{53} (2000), no. 2, 121--141.


\bibitem{DKhP} J. Ding, S. Khoroshkin, S. Pakuliak. Factorization of the
universal $R$-matrix for $U_q(\widehat{sl}_2)$ {\it Theoretical and
Mathematical Physics} {\bf 124:2} (2000), 1007-1036.


\bibitem{S} F. Smirnov. Form factors in completely integrable models of quantum
field theory, Adv. Series in Math. Phys., vol. 14, World Scientific, Singapore, 1992.


%\bibitem{T} V. Tarasov. A. Varchenko, Combinatorial Formulae for Nested
%Bethe Vectors.

\bibitem{VT} A. Varchenko, V. Tarasov. Jackson integrals representation of solutions to
the quantized
Knizhnik-Zamolodchikov equation, {\it St. Peterburg Math. Jour. } {\bf 6} (1995) no.2,
275--313.

\end{thebibliography}
\end{document}